\documentclass[11pt]{article}

\usepackage[top=1in,bottom=1.2in,left=1.25in,right=1.25in]{geometry}

\usepackage{graphicx}
\usepackage{ragged2e}
\usepackage{tabto}
\usepackage{amssymb}
\usepackage{float}
\usepackage{multirow}
\usepackage{rotating}
\usepackage{amsmath}
\usepackage[cmex10]{mathtools}
\usepackage{tikz}
\usepackage{amsthm}
\usepackage{color}
\usepackage{url}
\usepackage[noabbrev,capitalise]{cleveref}
\usepackage{afterpage}
\usepackage{multicol}
\usepackage[font=footnotesize]{caption}
\usepackage[font=footnotesize]{subcaption}
\usepackage{marginnote}
\usepackage{tcolorbox}
\usepackage{multirow}
\usepackage{newtxtext}
\usepackage{newtxmath}

\usepackage{multirow}

\usepackage[absolute,overlay]{textpos}

\usepackage[sf,bf,medium]{titlesec}

\usepackage{bm}
\usepackage{dsfont}
\usepackage{booktabs}
\usepackage{siunitx}
\usepackage{bm}

\newcommand\mtx[1]{\bm{\mathsf{#1}}}

\newcommand{\lp}{\left(}
\newcommand{\rp}{\right)}

\usepackage[sort,compress]{cite}

\newcommand\bx{\boldsymbol x}
\newcommand\by{\boldsymbol y}
\newcommand\bz{\boldsymbol z}
\newcommand\bE{\boldsymbol E}
\newcommand\bH{\boldsymbol H}
\newcommand\bA{\boldsymbol A}
\newcommand\bJ{\boldsymbol J}
\newcommand\bM{\boldsymbol M}
\newcommand\bN{\boldsymbol N}

\newcommand\bL{\boldsymbol L}
\newcommand\In{\operatorname{in}}
\newcommand\Sc{\operatorname{scat}}
\newcommand\bn{\boldsymbol n}
\newcommand\bu{\boldsymbol u}
\newcommand\bv{\boldsymbol v}

\newcommand\tot{\operatorname{tot}}
\newcommand\rhotail{\rho_{\textrm{tail}}}
\newcommand{\Ntri}{N_{\textrm{patches}}}
\newcommand{\nord}{n_{\textrm{order}}}
\newcommand{\niter}{n_{\textrm{iter}}}

\newtheorem{remark}{Remark}

\numberwithin{equation}{section}

\begin{document}


\begin{titlepage}

  \raggedleft
  {\sffamily \bfseries STATUS: arXiv pre-print}
  
  \hrulefill

  \raggedright
  \begin{textblock*}{\linewidth}(1.25in,2in) 
    {\LARGE \sffamily \bfseries Fast adaptive high-order integral equation
      methods\\
      for electromagnetic scattering from smooth\\ perfect electric conductors\\}
  \end{textblock*}

  \normalsize

  \vspace{2in}
  Felipe Vico\footnote{Research forms a part of the Advanced Materials program and was supported in part by MCIN with funding from the European Union NextGenerationEU (PRTR-C17.I1) and by Generalitat Valenciana, project: MAOCOM-6G, code: MFA/2022/056}\\
  \emph{\small Instituto de Telecomunicaciones y Aplicaciones
    Multimedia (ITEAM)\\
    Universitat Polit\`{e}cnica de Val\`{e}ncia\\
  Val\`{e}ncia, Spain}\\
  \texttt{\small felipe.vico@gmail.com}

   \vspace{\baselineskip}
Leslie Greengard\footnote{Research supported in part by
    the Office of Naval Research under award number~\#N00014-18-1-2307.}\\
  \emph{\small Courant Institute, NYU\\
    New York, NY 10012\\and\\Center for Computational Mathematics, Flatiron Institute\\
    New York, NY 10010.}\\
  \texttt{\small greengard@cims.nyu.edu}

   \vspace{\baselineskip}
  Michael O'Neil\footnote{Research supported in part by
    the Office of Naval Research under award
    numbers~\#N00014-17-1-2451, \#N00014-21-1-2383 and~\#N00014-18-1-2307.}\\
  \emph{\small Courant Institute, NYU\\
    New York, NY 10012}\\
  \texttt{\small oneil@cims.nyu.edu}
  
  \vspace{\baselineskip}
  Manas Rachh\\
  \emph{\small Department of Mathematics\\
    Indian Institute of Technology Bombay\\
  Mumbai, India}\\
  \texttt{\small mrachh@flatironinstitute.org}

  \begin{textblock*}{\linewidth}(1.25in,8.5in) 
    \today
  \end{textblock*}

\end{titlepage}

\begin{abstract}
  Many integral equation-based methods are available for problems of
  time-harmonic electromagnetic scattering from perfect electric conductors.
  Among the many challenges that arise in such calculations are the avoidance of
  spurious resonances, robustness of the method to scatterers of non-trivial
  topology or multiscale features, stability under mesh refinement, ease of
  implementation with high-order basis functions, and behavior in the static
  limit.  Since three-dimensional scattering is a challenging, large-scale
  problem, many of these issues have been historically difficult to
  investigate. It is only with the advent of fast algorithms for matrix-vector
  multiplies coupled with modern iterative methods that a careful study of these
  issues can be carried out effectively.  Our focus here is on comparing the
  behavior of several integral equation formulations with regard to the issues
  noted above, namely: the well-known, standard electric, magnetic, and combined
  field integral equations with standard RWG basis functions, and the more
  modern non-resonant charge-current and decoupled potential integral
  equation. Numerical results are provided to demonstrate the behavior of each
  of these schemes.  Furthermore, we provide some analytical properties and
  comparisons with the electric charge-current integral equation
  and the augmented regularized combined source integral equation.  \\
  

\end{abstract}

\small
\tableofcontents

\newpage

\normalsize
\section{Introduction}

Boundary integral equation methods are widely used in computational
electromagnetics (CEM), especially for exterior scattering problems. They impose
the outgoing radiation condition exactly, and for piecewise constant homogeneous
dielectrics or perfect conductors, reduce the dimensionality of the problem by
only requiring discretization of the boundary. While the number of degrees of
freedom required is dramatically reduced, these methods lead to dense linear
systems of equations -- hence, fast algorithms are needed to address large-scale
problems. At present, state-of-the-art solvers rely on iterative algorithms such
as GMRES or BiCGstab.  These algorithms work particularly well when the system
to be solved is well-conditioned with a spectrum that clusters away from the
origin.  When computing a solution via an iterative solver, each step requires a
matrix-vector product involving the system matrix. There are many algorithms
available for accelerating this step, and since it is by now fairly standard in
both academic and commercial software, we will rely here on the fast multipole
method (FMM)~\cite{MLFMM,Plane_Wave_FMM,HF-MLFMM,FMMBOOK,WIDEBAND}. Iterative
solvers with FMM acceleration only require an amount of work on the order
of~$\mathcal O(\niter N \log N)$ for any frequency, where $N$ is the system size
and~$\niter$ is the total number of iterations (and in the low-frequency regime,
the FMM solver generally scales as~$\mathcal O(\niter N)$).

The focus of the present paper is on the choice of integral formulation, its
discretization process, and the resulting effect on performance and accuracy
for scattering from closed surfaces. Different choices for each of
  these ingredients can sometimes result in dramatically difference results.  We
  will make a few brief remarks about open surfaces in the conclusions.
Currently, the most widely used solvers rely on the electric field, magnetic
field, and combined field integral equations (EFIE, MFIE and CFIE) discretized
using a Galkerin approach with RWG basis functions and a conforming mesh model
of the scatterer~\cite{RWG,MAUE,MautzHarrington}. These methods are subject to a
host of numerical difficulties, including low-frequency breakdown, high-density
mesh breakdown, and standard mathematical ill-conditioning. Rather than changing
the underlying formulation, the dominant approach in CEM has been to introduce
additional ideas to mitigate these problems. Loop-star basis
functions~\cite{Looptree2,Looptree3,ZC,andriulli2013well}, for example, improve
accuracy and conditioning in the low-frequency regime. Linear algebraic
pre-conditioners~\cite{Volakis_ILU,Andriulli_Calderon,helmholtz_stable_fast}
alleviate the difficulties produced by the hypersingular integral operator in
the EFIE, especially when dense meshes are needed to resolve sub-wavelength
features in the geometry. At the same time, there has been a significant effort
in the research community to develop {\em well-conditioned} Fredholm integral
equations of the second kind\footnote{Integral equations are said to be Fredholm
  equations of the second kind when the system matrix $A$ is of the form
  $I + K$, where $I$ is the identity operator and $K$ is a smoothing integral
  operator whose spectrum clusters at the origin. The condition number of such
  systems is typically independent of the number of degrees of freedom and
  stable under mesh refinement.}. While we do not seek to review the literature
here, these include the use of Calder\'on identities to \emph{analytically}
pre-condition the
EFIE~\cite{ChristiansenNedelec,Rokhlin_Maxwell,Andriulli_Calderon} and the use
of regularizing operators to pre-condition the
CFIE~\cite{CK2,Levadoux,Bruno_Krylow}. A complementary class of methods are the
so-called {\em charge-current formulations}. These methods are also aimed at
developing well-conditioned formulations that are free from low-frequency
breakdown~\cite{ZC,ZCCZ,CCIE1,CCIE3, CCIE4,CCIEour,Nystrom_our}, but achieve the
goal by introducing extra unknowns in the problem. Other formulations that lead
to resonance-free, second-kind equations include those based on generalized
Debye sources~\cite{gDEBYE,debye_boag} and decoupled potential
formulations~\cite{DPIE,Chew_A,DPIE_shanker}. More recently an augmented
regularized combined source integral equation (auRCSIE) was introduced
in~\cite{auRCSIE}. Rather than an exhaustive analysis of all such formulations,
we focus here on  two known representative integral equations from
  the existing literature (NRCCIE, and DPIE) comparing each with the
standard standard EFIE formulation discretized with RWG basis
functions~(EFIE-RWG).

Once an appropriate second-kind integral formulation has been selected, the
accuracy of the obtained solution will depend on the discretization and
quadrature methods used. In this paper, we investigate the use of a fast,
high-order, adaptive Nystr\"om-like method that yields high-order convergence
and permits adaptive refinement to capture small features in the geometry. Of
course, the quality of the geometry description itself also has an impact on the
accuracy of the results. Here, we use the method described
in~\cite{smooth_surface} which allows for the efficient construction of globally
smooth complex surfaces with multiscale features, high-order mesh generation,
and local refinement.

With all of this machinery in place, we are able to address challenging
electromagnetic scattering problems with millions of degrees of freedom in
physically delicate regimes. We show that for the effective solution of such
problems, {\em all} of the above ingredients play a role in robustly achieving
user-specified accuracies in the electric and magnetic fields: well-conditioned
formulations, high-order surface representations, and high-order quadratures
(complemented by suitable fast algorithms).

The key advances of this work are summarized below:
\begin{itemize}
\item A systematic comparison of existing integral formulations  for electromagnetic scattering from perfect conducting surfaces (see Table~\ref{tab1} for a summary). 
\item An efficient and high-order accurate FMM-accelerated implementation of DPIE, which to the best
of our knowledge is the only solver that avoids \emph{low-frequency breakdown},
  \emph{high-density mesh breakdown}, \emph{spurious resonances}, and is efficient even on 
  \emph{high-genus geometries}, and \emph{multiscale targets}.
  \item Demonstrations on a genus-17 geometry,
    a 16$\lambda$ cavity, and a 42$\lambda$ ship model
    confirm improved accuracy with fewer degrees of freedom.
\end{itemize}

\section{PEC integral equations with physical unknowns}

Electromagnetic scattering from a perfect electric conductor can be studied in
the time harmonic regime, where the full Maxwell equations reduce to:
\begin{equation} \nabla \times \bH^{\tot} = - i \omega\epsilon \bE^{\tot} \, , \quad 
 \nabla \times \bE^{\tot} =  i \omega \mu \bH^{\tot} .
\end{equation}
Here, we assume that the permittivity~$\epsilon$ and permeability~$\mu$ are
scalar constants. The perfect electric conductor (PEC) is defined by a bounded
region~$\Omega$
whose boundary is given by~$\Gamma = \partial \Omega$.
As is well-known, the boundary conditions on a PEC are~\cite{JACKSON,PAPAS}:
\begin{equation}
  \label{EHbc}
  \begin{aligned}
  \bn \times \bE^{\tot} &= {\boldsymbol 0}, &\qquad  \bn \cdot \bH^{\tot} &= 0, \\
  \bn \times \bH^{\tot} &= \bJ, &  \bn \cdot \bE^{\tot} &= \frac{\rho}{\epsilon}  ,
  \end{aligned}
\end{equation}
together with the continuity condition along the surface of the scatterer,
\begin{equation} i\omega \rho = \nabla_{\Gamma} \cdot \bJ .
\label{contcond}
\end{equation}
It is sufficient to enforce the boundary conditions on the tangential components
of the electric field, as done in the EFIE, but one or more of the other
(redundant) boundary conditions are often used in the alternative formulations
mentioned above and discussed below.

Furthermore, it is convenient to write the total field as the sum of a known
incoming field and an unknown scattered field:
\begin{equation}
  \bE^{\tot} = \bE^{\In} + \bE^{\Sc} \, , \quad 
  \bH^{\tot} = \bH^{\In} + \bH^{\Sc} \, . 
\end{equation}
The standard representation for the scattered fields is given in terms of a
vector and scalar potential,~$\mathbf{A,\phi}$, in the Lorenz gauge:
\begin{align}
   \bE^{\Sc} &= i \omega \bA^{\Sc} - \nabla \phi^{\Sc},  \label{Epotrep} \\
  \bH^{\Sc} &= \frac{1}{\mu} \nabla \times \bA^{\Sc}, \label{Hpotrep} 
\end{align}
with
\begin{equation}
\bA^{\Sc}[\bJ](\bx) =  \mu S_k[\bJ](\bx),
\end{equation}
\begin{equation}
 \phi^{\Sc}[\rho](\bx) = \frac{1}{\epsilon}  S_k[\rho](\bx),
\label{Skdef}
\end{equation}
and where $\bx\in\mathbb{R}^3\setminus \Omega$. The above layer potential
operators are defined by
\begin{equation}
  \begin{aligned}
  S_k[\mathbf{a}](\bx) &= \int_\Gamma g_k(\bx-\by) \, \mathbf{a}(\by) \, dA_{\by}, \\
  S_k[\sigma](\bx) &= \int_\Gamma g_k(\bx-\by) \, \sigma(\by) \, dA_{\by},
\end{aligned}
\end{equation}
with kernel given by the Green's function
\begin{equation}
g_k(\bx-\by) = \frac{e^{ik |\bx-\by|}}{ 4 \pi | \bx-\by|}  .
\end{equation}
Here, $\mathbf{a}$ is a tangential vector field and~$\sigma$ is a scalar-valued
function on the boundary $\Gamma$. It is important to note that the
charge~$\rho$ in~\eqref{Skdef} is not an extra degree of freedom, but must
satisfy the continuity condition~(\ref{contcond}). This ensures that the
resulting electromagnetic fields~$\bE^{\Sc},\bH^{\Sc}$ are Maxwellian.
Using the representation above for the scattered electric and magnetic fields,
the EFIE is obtained by imposing the boundary
condition~$\bn \times \bE^{\tot} = {\boldsymbol 0}$, the magnetic field integral
equation (MFIE) is obtained by imposing the boundary
condition~$\bn \times \bH^{\tot} = \bJ$, and the standard CFIE is obtained as a
linear combination of~$\bn \times \bH^{\tot} = \bJ$
and~$-\bn \times\bn \times \bE^{\tot} = {\boldsymbol 0}$. 

 Defining $\rho = \nabla_{\Gamma} \cdot \bJ/(i\omega)$ according to
  \eqref{contcond} causes the EFIE to be hypersingular and the evaluation of the
  electric field to be ill-conditioned.  In order to avoid these difficulties,
  charge-current formulations have been developed based on including electric
  charge as an independent unknown and imposing one (or more) additional
  conditions from~(\ref{EHbc}) in order to obtain a uniquely solvable system of
  equations~\cite{QianChew,CCIE1,CCIE3,CCIE4,CCIEour,Nystrom_our}. We will
often refer to this approach as an \emph{augmented formulation} since the number
of unknown and constraints have both increased. We turn now to the derivation of
two such scheme.

\subsection{Electric charge-current integral equation}

The electric charge-current integral equation (ECCIE) is presented
in~\cite{CCIEour}, following the ideas and nomenclature of~\cite{CCIE1, CCIE2}.
It is obtained from the representations~(\ref{Epotrep}) and~(\ref{Hpotrep}) by
imposing the conditions
\[ \bn \times \bH^{\tot} = \bJ, \quad
\bn \cdot \bE^{\tot} = \frac{\rho}{\epsilon} \, ,
\] 
yielding
\begin{equation}\label{NRCCIE_1}
 \frac{\bJ}{2}-M[\bJ]=\bn\times \bH^{\In}\\
 \end{equation}
and
\begin{equation}\label{NRCCIE_2}
  -i\omega\epsilon\mu\bn\cdot S_k[\bJ]+\frac{\rho}{2}+S'_k[\rho]=\epsilon\bn\cdot \bE^{\In} \, ,
\end{equation}
where 
\begin{equation}\label{Mdef}
 M[\bJ]= \bn \times \nabla \times S_k[\bJ] 
\end{equation}
is interpreted on surface in the principal value sense, and
\begin{equation}
  S'_k[\rho](\bx)=\int_{\Gamma}\frac{\partial g_k}{\partial n_{\bx}}(\bx-\by) \,
  \rho(\by) \, dA_{\by}.
\end{equation}
An analogous integral equation known as the Magnetic Charge-Current Integral
Equation (MCCIE) can be derived, but it shares similar properties and we will
not discuss the formulation in this paper.

\subsection{Non-resonant charge-current integral equation}

The non-resonant charge-current integral equation (NRCCIE) was introduced
in~\cite{CCIE1, CCIE2}, with a modified version in~\cite{CCIE4}. The basic idea
is to make use of~(\ref{NRCCIE_1}) and~(\ref{NRCCIE_2}), together with the
equation derived from imposing~${\bn \cdot \bE^{\tot} = \rho / \epsilon}$ and a
weak form of the continuity condition~(\ref{contcond}) obtained by integration
over the surface. These two equations take the form
\begin{equation}\label{NRCCIE_3}
  i\omega\mu\bn\times S_k[\bJ]-\frac{1}{\epsilon}\bn\times\nabla S_k[\rho]
  =-\bn\times \bE^{\In},
\end{equation}
\begin{equation}\label{NRCCIE_4}
  \nabla\cdot S_k[\bJ]-i\omega S_k[\rho]=0 \, .
\end{equation}
The NRCCIE is a system of two equations, the first obtained as a linear
combinations of~(\ref{NRCCIE_1}) and~(\ref{NRCCIE_2}), and the second obtained
as a linear combinations of~(\ref{NRCCIE_3}) and~(\ref{NRCCIE_4}):
\begin{equation}
    \label{NRCCIE}
    \frac{\bJ}{2}-M[\bJ] 
    + \alpha \bn\times\big\{ i\omega\mu \bn\times S_k[\bJ]-\frac{1}{\epsilon}\bn\times\nabla S_k[\rho]\big\} 
    = \bn\times \bH^{\In}-\alpha\bn\times\bn\times \bE^{\In}
\end{equation}
\begin{equation}
    \frac{\rho}{2}+S'_k[\rho]-i\omega\mu\epsilon \bn\cdot S_k[\bJ] +\alpha \big\{\nabla\cdot S_k[\bJ]-i\omega S_k[\rho]\big\}
= \epsilon\bn\cdot \bE^{\In} .
\end{equation}
Here, $\alpha$ is an arbitrary real positive constant. The NRCCIE is known to
have a unique solution at any frequency~\mbox{$\omega >0$}~(see
\cite{CCIE4}). The operators~$\nabla\cdot S_k[\bJ]$
and~$\bn\times\nabla S_k[\rho]$ are not compact, however, and therefore the
coupled system~(\ref{NRCCIE}) is not strictly speaking a Fredholm equation of
the second kind. Nevertheless, we will show that it has similar properties such
as a small condition number and the absence of high-density mesh breakdown.

\section{Decoupled potential integral equation}

Instead of solving for the physical quantities, current and charge, one can
instead indirectly solve for the vector and scalar potentials themselves. Such
an approach leads to the decoupled potential integral equation (DPIE),
originally introduced in~\cite{DPIE} to address the ubiquitous problem of
topological low-frequency breakdown endemic in almost all integral formulations
for electromagnetic scattering. The DPIE approach is based on considering two
uncoupled boundary value problems: one for the scalar potential, and one for the
vector potential. Trivially, both potentials satisfy the homogeneous Helmholtz
equation (due to the choice of Lorenz gauge).  For the scalar problem, consider
the boundary value problem:
\begin{equation}
    \begin{aligned}
\Delta \phi^{\Sc} +k^2 \phi^{\Sc} &=0\\
\phi^{\Sc}|_{\Gamma} - V &=-\phi^{\In}|_{\Gamma} \\
\int_{\Gamma} \frac{\partial \phi^{\Sc}}{\partial n} \, ds 
&=-\int_{\Gamma} \frac{\partial \phi^{\In}}{\partial n} \, ds,
    \end{aligned}
\end{equation}
where~$V$ is an unknown constant (voltage).   And
similarly, for the vector potential, consider the boundary value problem:
\begin{equation}
    \begin{aligned}
    \Delta \bA^{\Sc}+k^2\bA^{\Sc} &=0\\
    \bn\times \bA^{\Sc}|_{\Gamma}&=-\bn\times \bA^{\In}|_{\Gamma}\\
    \nabla \cdot \bA^{\Sc}|_{\Gamma} - v &=-\nabla \cdot 
    \bA^{\In}|_{\Gamma}\\
    \int_{\Gamma} \bn\cdot \bA^{\Sc} \, ds&=-\int_{\Gamma} \bn\cdot 
    \bA^{\In}\, ds,
    \end{aligned}
\end{equation}
where, as above,~$v$ is an unknown constant.  These boundary value problems
naturally extend to the case where~$\Omega$ is composed of multiple disjoint
components.  See~\cite{DPIE} for a thorough discussion of the role that the
constants~$V$ and~$v$
play in the representation of the fields. Each of these boundary value
problems can be solved by means of a second-kind integral equation using the
following representations for the scattered scalar and vector potentials:
\begin{equation}\label{fieldS1}
\phi^{\Sc}(\bx)=D_k[\sigma](\bx)-i\alpha S_k[\sigma](\bx),
\end{equation}
\begin{equation}
  \label{vecrep:DPIEv}
\bA^{\Sc}(\bx)=\nabla\times S_k[\mathbf{a}](\bx)- S_k[\bn\rho](\bx)\\
+i\alpha \big( S_k[\bn\times\mathbf{a}](\bx)+ \nabla S_k[\rho](\bx) \big),
\end{equation}
where we require that~$\alpha>0$ (but can be chosen freely), and where
\begin{equation}
  D_{k}[\sigma](\bx) =   \int_{\Gamma}\frac{\partial g_k}{\partial n_{\by}}(\bx-\by) \,
  \sigma(\by) \, dA_{\by}
\end{equation}
is the double layer potential.
Imposing the boundary conditions above, and using
the fact that~$\int_{\Gamma} D'_0 [\sigma] \, ds = 0$, see~\cite{DPIE},
eq. (A.11), we obtain the following system of equations for the
unknowns~$\sigma$,~$\boldsymbol{a}$,~$\rho$,~$V$, and~$v$:
\begin{equation}\label{HSS1}
\begin{aligned}
  \frac{\sigma}{2}+D_k[\sigma]-i\alpha S_k[\sigma]- V
    &=-\phi^{\In}|_{\Gamma}, \\
\int_{\Gamma} \big( (D'_k-D'_0)[\sigma]+i\alpha\frac{\sigma}{2} 
-i\alpha S'_k[\sigma] \big) ds 
&=-\int_{\Gamma} \frac{\partial \phi^{\In}}{\partial n}ds \, ,
\end{aligned}
\end{equation}
\begin{equation}
  \label{HVS1}
\frac{1}{2}
\left(\begin{array}{c}\mathbf{a} \\ \rho \end{array} \right) 
+
\mtx{L}
\left(\begin{array}{c}\mathbf{a} \\ \rho \end{array} \right) 
+
i\alpha \mtx{R} 
\left(\begin{array}{c}\mathbf{a} \\ \rho \end{array} \right) \\
+ \left(\begin{array}{c} 0 \\  v  \end{array} \right) 
=
\left(\begin{array}{c} -\bn\times \bA^{\In}|_{\Gamma} \\ 
  -\nabla \cdot \bA^{\In}|_{\Gamma} \end{array} \right),
\end{equation}
\begin{equation}
  \label{HVS2}
  \int_{\Gamma} \Big( -\bn\cdot S_k[\bn\rho]+i\alpha \big(\bn\cdot 
  S_k[\bn\times\mathbf{a}] \big)  \\
 -i\alpha\lp \frac{\rho}{2}+S'_k [\rho] \rp  \Big) ds 
=-\int_{\Gamma} \bn \cdot \bA^{\In} \, ds \, ,
\end{equation}
The matrix integral operators $\mtx{L}$ and~$\mtx{R}$ above are defined by:
\begin{equation}
\begin{aligned}
\mtx{L} 
\left(\begin{array}{c}\mathbf{a} \\ \rho \end{array} \right) 
=&\left( \begin{array}{c}
L_{11}[\mathbf{a}] + L_{12}[\rho] \\
L_{21}[\mathbf{a}] + L_{22}[\rho]
\end{array} \right),
\quad 
\end{aligned}
\end{equation}
where
\begin{equation}\label{operators_vect}
\begin{aligned}
L_{11}[\mathbf{a}]=&\hat{\bn} \times S_k [\mathbf{a}], \\
L_{12}[\rho]=&\ -\hat{\bn}\times S_k[\hat{\bn}\rho]),\\
L_{21}[\mathbf{a}]=&\ 0,\\
L_{22}[\rho]=&\ D_k[\rho],\\
\end{aligned}
\end{equation}
and
\begin{equation}
\begin{aligned}
\mtx{R}
\left(\begin{array}{c}\mathbf{a} \\ \rho \end{array} \right) 
=&\left( \begin{array}{c}
R_{11}[\mathbf{a}] + R_{12}[\rho] \\
R_{21}[\mathbf{a}] + R_{22}[\rho]
\end{array} \right),
\end{aligned}
\end{equation}
where
\begin{equation}
\begin{aligned}
R_{11}[\mathbf{a}]=&\ \hat{\bn}\times S_k[\hat{\bn}\times\mathbf{a}],\\
R_{12}[\rho]=&\ \hat{\bn}\times\nabla S_k[\rho],\\
R_{21}[\mathbf{a}]=&\ \nabla \cdot S_k[\hat{\bn}\times\mathbf{a}],\\
R_{22}[\rho]=&\ -k^2S_k[\rho].
\end{aligned}
\end{equation}

The vector integral equation above in~\eqref{HVS1} and~\eqref{HVS2} is not,
strictly speaking, a Fredholm equation of the second kind since~$R_{12}$
and~$R_{21}$ are bounded but \emph{not} compact operators. Nevertheless, we will
show that it has similar properties. The formulation is resonance free and
stable at arbitrarily low frequencies for geometries of any genus having
multiple components (see~\cite{DPIE} for further detail). The original
formulation in~\cite{DPIE} contains an additional regularizing operator that we
have omitted here for simplicity. Stability does not appear to be compromised in
our experiments. The coefficient~$\alpha$ is included above to avoid spurious
resonances; we typically set~$\alpha=1$, but for complicated geometries, it may
be possible to optimize the choice in order to reduce the total number of
iterations.
\begin{remark}
If~$\alpha=0$, we will refer to the resulting (simpler) integral equation as the
resonant DPIE (rDPIE). The spurious resonances are actually the same as those
for the MFIE.
\end{remark}

\section{Properties of various integral formulations}

We summarize the expected properties (based on a mathematical analysis) of
the various formulations in the table below.
We further describe some of the items in the left-hand column of Table~\ref{tab1}:
\begin{itemize}
  \item A \emph{spurious resonance} is a frequency where the integral equation
        is not invertible but the scattering problem is itself well-posed.
  \item \emph{High-density mesh breakdown} refers to a significant growth in the
        numerical condition number of the finite-dimensional linear system to be
        solved under mesh refinement. Some integral equations are Fredholm
        equations of the second kind which, in the absence of spurious
        resonances, have bounded condition numbers independent of the number of
        degrees of freedom.
  \item \emph{Catastrophic cancellation} in~$\bE^{\Sc}$,~$\bH^{\Sc}$ refers to a
        loss of precision in computing the scattered fields of interest once the
        integral equation has been solved (see section~\ref{farfield-sec}).
  \item Second kind integral equations and equations whose system matrices are
        of the form $I+K$, where $K$ is the discretiztion of a bounded operator,
        tend to converge rapidly using GMRES or BiCGSTAB as an iterative method.
  \item An equation is \emph{stable at low frequency} if the condition number
        does not grow as the frequency tends to zero. This can be the case for
        surfaces without holes (of genus zero) or more generally (for surfaces
        of arbitrary genus).
\end{itemize}

  \begin{table}[!b] 
    \centering
    \caption{Properties of various integral equation formulations}
    \begin{tabular}{c|cccccccc} 
      &
    \begin{turn}{90}EFIE\end{turn} &
    \begin{turn}{90}MFIE\end{turn} &
    \begin{turn}{90}CFIE\end{turn} &
    \begin{turn}{90}ECCIE\end{turn} &
    \begin{turn}{90}NRCCIE\end{turn} &
    \begin{turn}{90}DPIE\end{turn}\\
    \hline
    \begin{tabular}{@{}c@{}}Resonance-free \end{tabular} 
      &   &   & \checkmark &   &   & \checkmark \\
    \hline
    \begin{tabular}{@{}c@{}} No high-density \\ mesh breakdown\end{tabular}
      &   & \checkmark &   & \checkmark & \checkmark & \checkmark \\
    \hline
    \begin{tabular}{@{}c@{}} Free from catastrophic\\cancellation in
    $\bE^{\Sc}$, $\bH^{\Sc}$ \end{tabular} 
      &   &   &   & \checkmark & \checkmark & \checkmark \\
    \hline
    \begin{tabular}{@{}c@{}} Second kind \\ Fredholm eq. \end{tabular} 
      &   & \checkmark &   & \checkmark &   & \checkmark \\
    \hline
    \begin{tabular}{@{}c@{}} Rapid convergence with \\ GMRES, BiCGSTAB
    \end{tabular} 
      &   & \checkmark &   & \checkmark & \checkmark  & \checkmark \\
    \hline
    \begin{tabular}{@{}c@{}} Stable at low frequencies \\
    for surfaces of genus zero \end{tabular} 
      &   & \checkmark &   & \checkmark & \checkmark  & \checkmark \\
    \hline
    \begin{tabular}{@{}c@{}} Stable at low frequencies \\ for surfaces 
    of any genus \end{tabular}
      & & & & & & \checkmark \\
    \hline
    \end{tabular}
    \label{tab1}
    \end{table}
    
    Following a description of discretization schemes, subsequent sections of the paper
provide numerical evidence that the properties summarized above have practical consequences.

\section{Surface representation, discretization, and quadrature}
\label{intandinterp}

Given a surface that has been approximated using flat triangles,
it is a standard procedure to discretize the EFIE, MFIE
or CFIE using edge-based RWG basis functions~\cite{RWG} in a Galerkin framework;
this corresponds to linear current profiles on each triangle. Since the
formulation is standard, we will not describe it in further detail. We will also
investigate the performance of higher-order non-Galerkin discretizations. In
this case, we must also assume that the surface~$\Gamma$ is described as a set
of triangular patches~\mbox{$\Gamma=\cup_{j=1}^{\Ntri} \Gamma_j$}, where $\Ntri$
is the number of curved triangular patches. For each $j$, we assume there exists
a known parameterization~$\bx_{j}$ such that
\begin{equation}
    \bx_j :T\rightarrow \Gamma_j\subset\mathbb{R}^3,
\end{equation}
where~$T$ is the canonical unit triangle:
\begin{equation}
\begin{aligned}
    T=\big\{(u,v) : u\ge0,v\ge 0,u+v \le 1\big\} \, .
\end{aligned}    
\end{equation}

\begin{remark}
  Since many computer-aided design systems or meshing algorithms produce only
  flat triangulations, the surfaces used as examples in this paper are generated
  using the algorithm of~\cite{smooth_surface}. This results in a surface of the
  desired form, with the regularity (curvature) of the surface locally
  controlled, permitting adaptive refinement and resolution of multiscale
  features.
\end{remark}

Given the surface $\Gamma$ described by an atlas of functions
$\left\{ \bx_j \right\}$, we also require a suitable set of sampling/quadrature
nodes and a suitable representation of smooth functions on each~$\Gamma_j$. For
this task, we will use Vioreanu-Rokhlin nodes/weights~\cite{vioreanu_rokhlin}
and Koornwinder polynomials as a basis for smooth functions, respectively. The
Vioreanu-Rokhlin nodes and weights have been designed so that the quadrature
rule
\begin{equation}
  \label{numquad}
    \int_T f(u,v) \, du \, dv \approx\sum_{i=1}^p w_{i} \, f(u_i,v_i),
\end{equation}
with 
\begin{equation}
  p=(\nord+1)(\nord+2)/2
\end{equation}
nodes exactly integrates (to machine precision) all polynomials~$u^k v^l$ in two
variables with total degree satisfying~\mbox{$k+l\le \nord$}. Furthermore, the
Koornwinder polynomials~$P_{n,m}$ on the standard triangle~$T$ are given
explicitly by
\begin{equation}
P_{n,m}(u,v)=(1-v)^m  P^{0,2m+1}_{n-m}(1-2v)\, \, P^{0,0}_m\left(\frac{2u}{1-v}-1\right),
\end{equation}
with $n=0,1,2,...$ and $m = 0, 1, \ldots, n$. Here $P_n^{\alpha,\beta}$ for
${n\in \mathbb{N}}$ are the standard Jacobi polynomials which are orthogonal with
respect to the weight function~$(1-x)^{\alpha}(1+x)^{\beta}$ on the interval
[-1,1], see~\cite{nist}. This is an orthogonal basis that comes equipped with
fast and stable recurrence formulas for their evaluation. Moreover, the mapping
from samples of functions at Vioreanu-Rokhlin nodes to the corresponding
coefficients in the Koornwinder basis is well-conditioned and straightforward to
generate. We refer the reader to~\cite{vioreanu_rokhlin,koornwinder,fmm-quad}
for further details. Once the Koornwinder expansion of a function is available,
it is a simple matter of evaluation to interpolate that function with high-order
accuracy to any other point on the triangle.

Additionally, we also require a basis in which to describe tangential vector
fields along each patch. To this end, we construct two sets of vector-valued
basis functions on each patch $\Gamma_{j}$ as follows. We first set
\begin{equation}
\begin{aligned}
    \bu_j(u,v)&=\frac{\partial \bx_j}{\partial u} \\
    \bn_j(u,v)&=\frac{\partial \bx_j}{\partial u}\times \frac{\partial \bx_j}{\partial v} \\
\end{aligned}    
\end{equation}
and 
\begin{equation}
\begin{aligned}
    \hat{\bu}_j(u,v)&=\frac{\bu_j(u,v)}{|\bu_j(u,v)|} \\
    \hat{\bn}_j(u,v)&=\frac{\bn_j(u,v)}{|\bn_j(u,v)|} \\
    \hat{\bv}_j(u,v)&=\hat{\bn}_j(u,v)\times\hat{\bu}_j(u,v).
\end{aligned}    
\end{equation}
Clearly~$\hat{\bu}$, $\hat{\bv}$, $\hat{\bn}$, form a pointwise orthonormal set
of coordinates along~$\Gamma_{j}$. Then, we set
\begin{equation}
\begin{aligned}
    \mathbf{U}^j_{n,m}(u,v) = P_{n,m}(u,v) \, \hat{\bu}_j(u,v)\\
    \mathbf{V}^j_{n,m}(u,v) = P_{n,m}(u,v) \, \hat{\bv}_j(u,v).
\end{aligned}    
\end{equation}
These vector basis functions are furthermore orthonormal in the sense that
\begin{equation}
\begin{aligned}
  \int_{T}\mathbf{U}^j_{n,m}(u,v)\cdot\mathbf{U}^j_{n',m'}(u,v) \, du \, dv& = \delta_{n,n'}\delta_{m,m'} ,\\ \int_{T}\mathbf{V}^j_{n,m}(u,v)\cdot\mathbf{V}^j_{n',m'}(u,v) \, du \, dv &=\delta_{n,n'} \delta_{m,m'} ,\\
  \int_{T}\mathbf{U}^j_{n,m}(u,v)\cdot\mathbf{V}^j_{n',m'}(u,v) \, du \, dv&=0  .
\end{aligned}    
\end{equation}
In our method, tangential vector fields are represented at each Vioreanu-Rokhlin
node~$(u_i,v_i)$ via an expansion in the two sets of basis functions
$\{ \mathbf{U}^j_{m,n}, \mathbf{V}^j_{m,n} \}$ (and evaluated at additional
points on~$\Gamma_j$ as needed using the Koornwinder basis).

\begin{remark}
  The reader may have noted that the basis functions used to discretize
  tangential vector fields, such as the electric current, do not correspond to
  a \emph{div-conforming} discretization~. Indeed, no continuity of any kind is
  enforced between adjacent triangles. This makes discretization very
  straightforward, as it can be done independently for each triangular patch. As
  we will see in the numerical examples, this choice does not introduce any
  artifacts, even at second-order accuracy. The robustness of the method is due
  to the accuracy of the integration method described below and a fundamental
  fact about second-kind integral equations: when using Nystr\"om
  discretizations, \emph{the order of accuracy of the method is equal to the
    order of accuracy of the underlying quadrature
    scheme}~\cite{anselone_book,Kress,ATKINSON}. A Nystr\"om method is one in which the
  integral equation is converted to a finite dimensional linear system by merely
  sampling the kernel and the unknown at a collection of quadrature nodes and
  approximating the integral operator by a quadrature rule over those same
  nodes~\cite{ATKINSON}.
\end{remark}

\subsection{Near and far field quadrature}

Since the integral operators appearing in all of our representations are
non-local, it is convenient to make use of a quadrature scheme that exploits the
smoothness of the integrand in various disjoint situations. If
the integral is taken over some triangle~$j$ (which we will call the {\em
  source} triangle) and the target point is located on triangle~$i$, we have the
following separate regimes: Following the discussion in~\cite{fmm-quad}, we
distinguish the {\em self-interaction} (when $i=j$), the near field (when
$i\neq j$ but the triangles are adjacent or nearby), and the far field (when
$i\neq j$ and the triangles are far apart). For the far field interactions, we
use the Vioreanu-Rokhlin quadrature described above (suitable for smooth
functions) and for which the fast multipole method FMM \cite{MLFMM},
\cite{Plane_Wave_FMM} can be applied directly to the discrete sum to accelerate
the computation. The self interaction is computed by a specialized high-order
quadrature rule due to Bremer and Gimbutas \cite{Singular_spq_2}.

The near interactions correspond to an integrand which is formally smooth but
very sharply peaked at the target. These are, in some sense, the most cumbersome
integrals to evaluate. For these, we rely on the method introduced in
\cite{fmm-quad}, which uses adaptive quadrature on the source
triangle~$\Gamma_j$ with a carefully precomputed multiscale hierarchy of
interpolants for the underlying density to reduce the cost. (Although the cost
is linear in the total number of degrees of freedom, the accurate evaluation of
near field quadratures is the most expensive step in quadrature generation.)

\begin{remark}
  In the end, each integral operator such as $S_k[\sigma](\bx)$ in
  \eqref{Skdef}, can be approximated in the form of a matrix with entries
  \[
    \begin{aligned}
      S_k^{disc}(i,j) &= g_k(\bx_i - \bx_j) w_{ij}, &\qquad &{\rm for\ } i\neq
                                                              j, \\
      S_k^{disc}(i,i) &= w_{ii} , & &
    \end{aligned}
  \]
  which maps point values of $\sigma$ at the discretization nodes to the value
  of the integral at those same nodes. The formulae for~$w_{ij}$ can be rather
  involved when points~$i$ and~$j$ lie on the same triangle or nearby-triangles,
  and we refer the reader to the literature cited above for details.
\end{remark}

\subsection{Error estimation}
\label{sec:err-est}
The use of orthogonal basis functions to represent the source densities on each
triangle has an additional advantage beyond high-order accuracy itself. Namely,
these representations can be used for \emph{a posteriori} error estimation and
as a monitor for identifying regions which need further geometric refinement.
The procedure is straightforward: from the samples of the unknown densities on
each patch, we obtain the coefficients of the corresponding function
approximation in the Koornwinder basis. This basis has the property that a
well-resolved function has a rapid decay of its Koornwinder coefficients. A
basic heuristic for the local error is to simply examine the relative norm of
the highest order basis functions. More precisely, let us first consider a
scalar quantity, such as the induced charge~$\rho$ on~$\Gamma_j$. From the
discussion above, using the Nystr\"om-like method, after solving our integral
equation we have the discrete values~$\rho(\bx_j(u_i,v_i))=\rho_{ji}$ at the
Vioreanu-Rokhlin nodes. Let us denote the corresponding Koornwinder
approximation by:
\begin{equation}
  \rho(\bx)=\rho(\bx_j(u,v))\approx\sum_{m+n\le \nord} c^j_{n,m} \, P_{n,m}(u,v).
\end{equation}
We may then define the following function as our error monitor on this patch:
\begin{equation}
  \begin{aligned}
    \rhotail(\bx) &=\rhotail(\bx_j(u,v)) \\
    &=\sum_{m+n= \nord} c^j_{n,m} \, P_{n,m}(u,v)
  \end{aligned}
\end{equation}
with $L^2$ norm
\begin{equation}
    \delta\rho_j=\sqrt{\int_{\Gamma_j}|\rhotail(\bx)|^2 \, dA_{\bx}}
\end{equation}
which serves as a triangle-by-triangle error estimate. The global absolute and
relative errors can then be estimated as
\begin{equation}
    \| \delta \rho \|_2\approx\sqrt{\int_{\Gamma}|\rhotail(\bx)|^2 dA_{\bx}}   
\end{equation}
and 
\begin{equation}
  \frac{\|\rho_{\text{error}}\|_2}{\|\rho\|_2}
  \approx\sqrt{\frac{\int_{\Gamma}|\rhotail(\bx)|^2 \,
      dA_{\bx}}{\int_{\Gamma}|\rho(\bx)|^2 \, dA_{\bx}}} \, ,
\end{equation}
respectively. 
Such estimates are uniformly robust for second-kind integral 
equation formulations \cite{boyd,lee1997fast,trefspec}.
As we will see below, the far field errors are approximately one
order of magnitude smaller. This is not
surprising, since the field quantity is obtained from the density through the
process of integration. Note that the $\ell_2$ norm of the sequence
$\{\delta\rho_j\}_{j=1}^{\Ntri}$ equals $\|\delta \rho\|_2$. Plotting the
piecewise constant function $\delta\rho_j$ on the triangulated surface helps
visualize regions with large errors and identifies triangles which require local
refinement if the obtained accuracy is not sufficient. The error estimation is
analogous for vector densities, such as the electric current.

\section{Far field estimation}
\label{farfield-sec}

The far field induced by a given electric or magnetic current can be computed
from the Fourier transform of the currents themselves (see \cite{cardama}). In
some of our formulations, such as the DPIE, the unknowns are
non-physical quantities. One could develop expressions for the far field in
terms of these unknowns using standard parallel-ray approximations. This
approach, however, has some disadvantages that we will discuss later. A second
option is to use a spherical \emph{proxy surface} that contains the full
scatterer and first compute the corresponding electric and magnetic fields on
that sphere. The principle of equivalent currents can then be used to compute
the field at any point in the far field (or the far field pattern itself). This
latter method has some stability advantages, and is worth describing in more
detail.

For known electric and magnetic currents~$\bJ,\bM$ along the proxy
sphere~$S_{R_0}$ of radius~$R_0$ (with~$R_0$ sufficiently large so as to enclose
the scatterer), the far field pattern is given by:
\begin{equation}
    \begin{aligned}
    E_{\theta}(\hat{\bx})&=i\frac{e^{ik|\bx|}}{2\lambda |\bx| }\Big(\eta N_{\theta}(\hat{\bx})+L_{\phi}(\hat{\bx})\Big),\\
    E_{\phi}(\hat{\bx})&=i\frac{e^{ik|\bx|}}{2\lambda |\bx| }\Big(\eta N_{\phi}(\hat{\bx})-L_{\theta}(\hat{\bx})\Big),\\
    H_{\phi}(\hat{\bx})&=\frac{1}{\eta}E_{\theta}(\hat{\bx}),\\
    H_{\theta}(\hat{\bx})&=-\frac{1}{\eta}E_{\phi}(\hat{\bx}),\\
    \end{aligned}
\end{equation}
where
\begin{equation}
  \label{Fourier_NF}
  \begin{aligned}
    \bN(\hat{\bx})&=\int_{S_{R_0}}\bJ(\by) \, e^{-ik\hat{\bx}\cdot\by} \,
                    dA_{\by}, \\
    \bL(\hat{\bx})&=\int_{S_{R_0}}\bM(\by) \, e^{-ik\hat{\bx}\cdot\by} \,
                    dA_{\by}  .
  \end{aligned}
\end{equation}
Above,~$\lambda$ denotes the wavelength~$\lambda = 2\pi/k$ and~$\eta =
\sqrt{\mu/\epsilon}$, the free-space impedance.
  The relevant currents can be computed on~$S_{R_0}$ from the scattered
fields~$\bE, \bH$ (using the FMM for efficiency) according to the equivalent
current principle:
\begin{equation}
    \begin{aligned}
        \bJ&=\hat{\bn}\times\bH , \\
        \bM&=-\hat{\bn}\times\bE  .
    \end{aligned}
\end{equation}
Here~$\hat{\bn}$ is the outward unit normal to the sphere~$S_{R_0}$. If the
scatterer is electrically large, the projection integrals in (\ref{Fourier_NF})
are expensive to evaluate naively by direct quadrature over a sufficiently fine
mesh on $S_{R_0}$. In that case, the fast Fourier transform (FFT) or its
non-uniform variant (NUFFT) can be used to accelerate the calculation
\cite{boag_fast_po,NUFFT1,NUFFT2}.

Unfortunately, the expressions in~(\ref{Fourier_NF}) are unstable at
low-frequency and subject to catastrophic cancellation. This problem is
discussed in~\cite{CCIEour} and stems from the fact that the magnitude of the
far field is $O(\omega)$ while the integrand is $O(1)$. The stabilization
introduced in~\cite{CCIEour} is based on introducing equivalent electric and
magnetic charges. These equivalent charges can easily be obtained from the
normal components of the fields~$\bE$,~$\bH$ on the spherical proxy surface.
Numerically stable (and exact) expressions for~$\bN$ and~$\bL$ are then given by
\begin{equation}
  \label{Fourier_NF_lf}
  \begin{aligned}
    \bN(\hat{\bx}) &= \int_{S_{R_0}} \lp \bJ(\by) (e^{-ik\hat{\bx}\cdot\by}-1) 
      -i\omega \, \by \, \rho(\by) \rp  dA_{\by}, \\
    \bL(\hat{\bx})&=\int_{S_{R_0}} \lp \bM(\by) (e^{-ik\hat{\bx}\cdot\by}-1) -i\omega \, 
    \by \, \sigma(\by) \rp  dA_{\by},
    \end{aligned}
\end{equation}
where
\begin{equation}
    \begin{aligned}
        \sigma&=\hat{\bn}\cdot\bH, \\
        \rho&=\hat{\bn}\cdot\bE  .
    \end{aligned}
\end{equation}
Note that the term $(e^{-ik\hat{\bx}\cdot\by}-1)$ is also of the order
$O(\omega)$ and can be evaluated without catastrophic cancellation as
\begin{equation}
  (e^{-ik\hat{\bx}\cdot\by}-1)
  =2ie^{i\frac{k}{2}\hat{\bx}\cdot\by}
  \sin\left( \frac{k}{2}\hat{\bx}\cdot\by
  \right) .
\end{equation}
In short, the expressions in~\eqref{Fourier_NF} are slightly more accurate at high
frequencies, while the expressions in~\eqref{Fourier_NF_lf} are significantly more
accurate and stable at low frequencies. Thus, we recommend the use
of~\eqref{Fourier_NF} for scatterers that are larger than $0.5$ wavelengths in
size and~\eqref{Fourier_NF_lf} otherwise.

\section{Numerical examples}

\begin{figure}[!t]
  \centering
  \includegraphics[width=.6\linewidth]{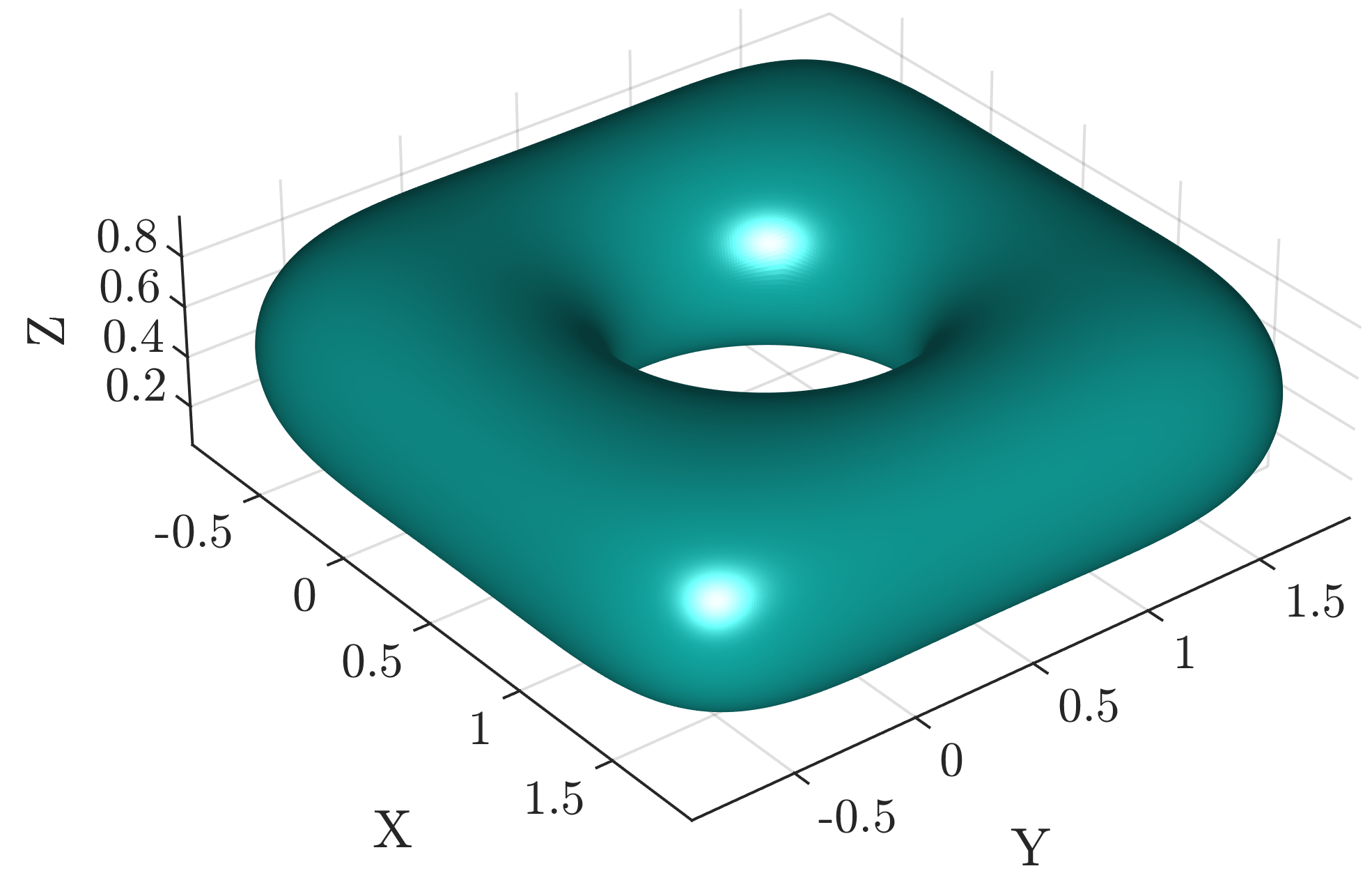}%
  \caption{A smoothed rectangular torus of genus one.}
  \label{Toroidal_Surface}
\end{figure}

In this section, we illustrate the behavior of the integral representations and
discretization methods discussed in the preceding sections. For
sections~\ref{subsec:convergence}, and~\ref{subsec:iterative}, the scatterer is
either a sphere of radius~$R=1$m or a smooth version of a rectangular torus, see
Fig.~\ref{Toroidal_Surface}. The toroidal geometry was obtained via the
surface smoothing algorithm of~\cite{smooth_surface} applied to a rectangular
torus defined as the union of rectangular faces parallel to the coordinate axes.

The code was implemented in Fortran and compiled using the GNU Fortran
11.2.0 compiler. We use the fast multipole method implementation from
the~\texttt{FMM3D}
package\footnote{\url{https://github.com/flatironinstitute/FMM3D}}, the high
order local quadrature corrections from the~\texttt{fmm3dbie}
package\footnote{\url{https://github.com/fastalgorithms/fmm3dbie}}, and the
high-order mesh generation code from the~\texttt{surface-smoother}
package\footnote{\url{https://github.com/fastalgorithms/surface-smoother}}.

In each of the examples, unless stated otherwise, the surface is represented
using flat triangles and the integral equations are discretized using a Galerkin
approach with the Rao-Wilton-Glisson basis and test
functions~\cite{RWG,MFIE-RWG} for the EFIE, MFIE, and CFIE. On the other hand,
for the NRCCIE and DPIE, the surface is represented using a collection of
high-order curvilinear triangles, and the integral equations are discretized
using a Nystr\"{o}m-like approach with locally-corrected quadratures.

\subsection{Accuracy}
\label{subsec:convergence}

To test the accuracy of the solvers, we can generate nontrivial exact solutions to the
boundary value problem (i.e. the scattering problem) and validate our numerical
approximations. For this, we define the true solution as the electromagnetic field induced by
a magnetic dipole (or set of magnetic dipoles) in the interior of the scatterer.
This induces boundary data which is then provided to the integral equation solver.
Evaluating the computed solution can then be compared to the kn own dipole field.
We define ~$\varepsilon_{E}$
and~$\varepsilon_{H}$ as the relative $L^2$ error in the electric and
magnetic fields at a collection of targets in the exterior region, and
let~$\varepsilon_{a}= \max(\varepsilon_{E},\varepsilon_{H})$. When the conductor
is a sphere, we may also check the accuracy of the computed scattered fields generated by
an incident plane wave, since an exact solution in the exterior is given by the
Mie series. With a slight abuse of notation, we will use~$\varepsilon_{a}$ to
denote this error as well.

\subsubsection{Convergence}

\begin{figure}[!t]
\begin{center}
\includegraphics[width=0.8\linewidth]{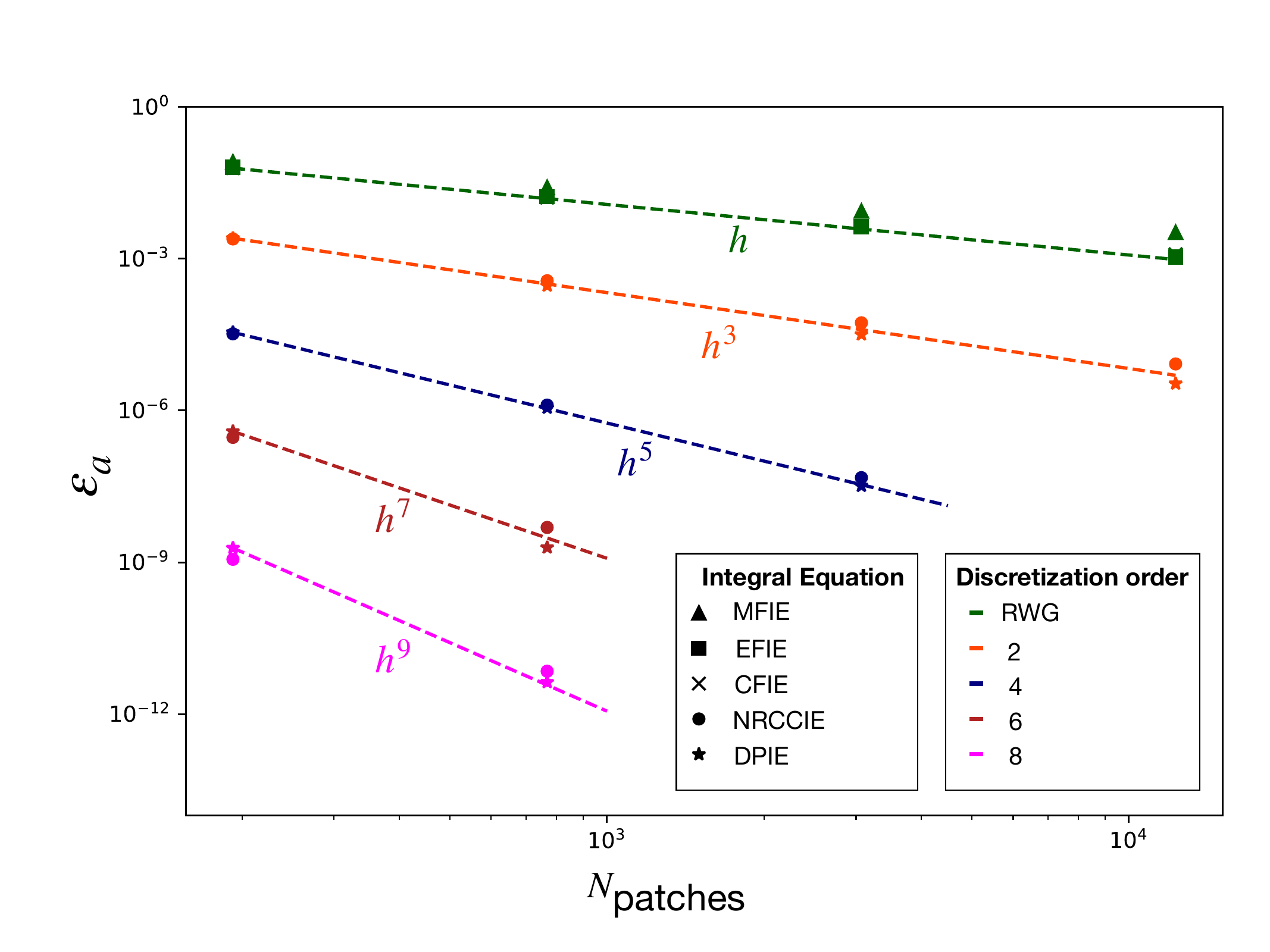}%
\end{center}
\caption{Relative error in the scattered field of a p.e.c. sphere of diameter $ D=\frac{1}{\pi}\lambda$ and incoming linearly polarized plane wave. We compare NRCCIE, integral equations with discretization order 2,4,6 and, 8 and the standard CFIE, MFIE, EFIE integral equations discretized with RWG basis functions}
\label{fig:sphere_convergence}
\end{figure}

In Fig.~\ref{fig:sphere_convergence}, we plot the error $\varepsilon_{a}$
corresponding to scattering from a PEC sphere with radius~$R=\SI{1}{m}$ and
wavenumber~$k=\SI{1}{m^{-1}}$ (the diameter of the sphere is $\lambda/\pi$) due
to an incoming linearly polarized planewave for each of the EFIE, MFIE, CFIE,
NRCCIE, and DPIE; results for the EFIE, MFIE, and CFIE are reported using RWG
basis functions, and results for the NRCCIE and DPIE are reported for
discretization orders $\nord = 2,4,6,8$. The errors decrease at the expected
rate of~$O(h)$ for the EFIE, MFIE, and CFIE, and at the expected rate
of~$O(h^{\nord+1})$ for NRCCIE and DPIE. Here $h$ is the diameter of a typical
triangle in the discretization.

\subsubsection{Absence of spurious resonances}
\label{subsec:spurious}

\begin{figure}[!t]
  \centering
  \includegraphics[width=.8\linewidth]{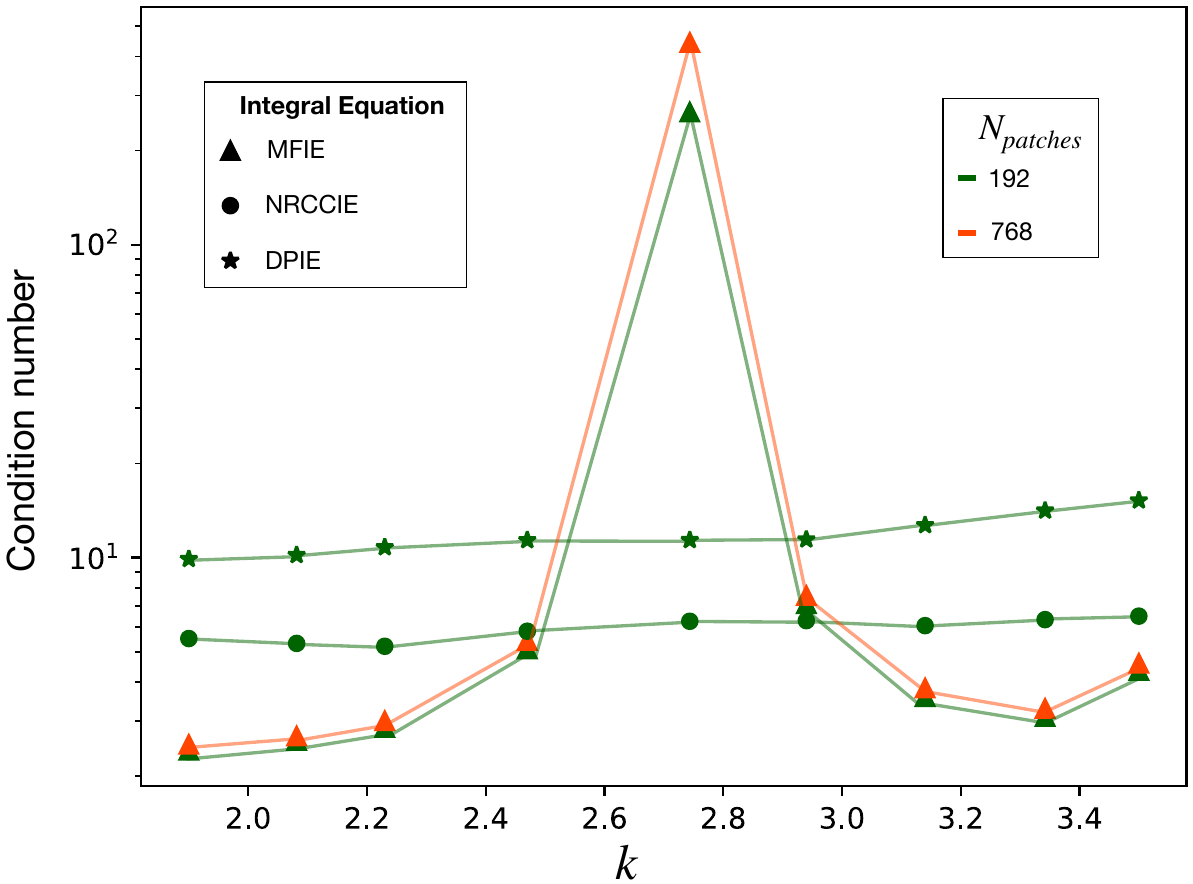}
  \caption{Condition number of the discretized integral equations for the MFIE, NRCCIE, and DPIE.}
  \label{fig:spurious}
\end{figure}

The exterior scattering problem has a unique solution for any real
wavenumber~$k$. However, it is well-known that the MFIE has spurious resonances,
i.e. wavenumbers~$k$ for which the integral equation is not invertible. On the
sphere, these spurious resonances can be computed analytically. To demonstrate
the absence of spurious resonances for the NRCCIE, and DPIE, we plot
the condition number of the discretized integral equations as a function of~$k$
in Fig.~\ref{fig:spurious}. All of the integral equations were discretized
using 192 patches and $\nord=2$. 
The interval~$k\in [1.9,3.5]$ has one internal
resonance of the MFIE on the sphere of radius $R=\SI{1}{m}$ given by
\begin{equation}
  \begin{aligned}
    k_{1} &= \SI{2.743707269992265}{m^{-1}} \, . \\
  \end{aligned}
\end{equation}
We observe that all integral equations except the MFIE have a bounded condition
number on the range of values of~$k$ considered, while the MFIE has a high
condition number precisely at its spurious resonant wavenumber. To further
confirm the presence of the spurious resonance, we also plot the condition
number for the MFIE using $768$ patches and observe that the condition number of
the resulting system increases as we obtain a more accurate discretization of
the integral equation at the spurious resonance, while there is very little
impact on the condition number at the other wavenumbers.  When computing the
condition numbers of the discretized linear systems, we scale both the unknowns
and the boundary data using the square root of the smooth quadrature weights to
obtain a better approximation of the integral equation in an $L^{2}$
sense~\cite{bremer_zydrunas}.

\subsubsection{Static limit}
\label{subsec:num-staticlimit}

In the static limit, the boundary value problems for the electric and magnetic
fields completely decouple. The fields computed at finite, but small
wavenumbers, converge to the solutions of the boundary value problems for the
electrostatic and magnetostatic fields. Since there exists a stable limit for
the underlying system of differential equations, it is a desirable feature that
the integral equations remain stable in the static limit as well. Integral
equation methods tend to have two kinds of failure modes in the static limit:
(1) deterioration in the accuracy of the computed solution using a fixed
discretization which resolves both the geometry and the boundary data
as~$k \to 0$; and (2) failure to converge at the expected rate upon mesh
refinement for a fixed, but small~$k$.

\begin{figure}[t]
\begin{center}
  \includegraphics[width=0.8\linewidth]{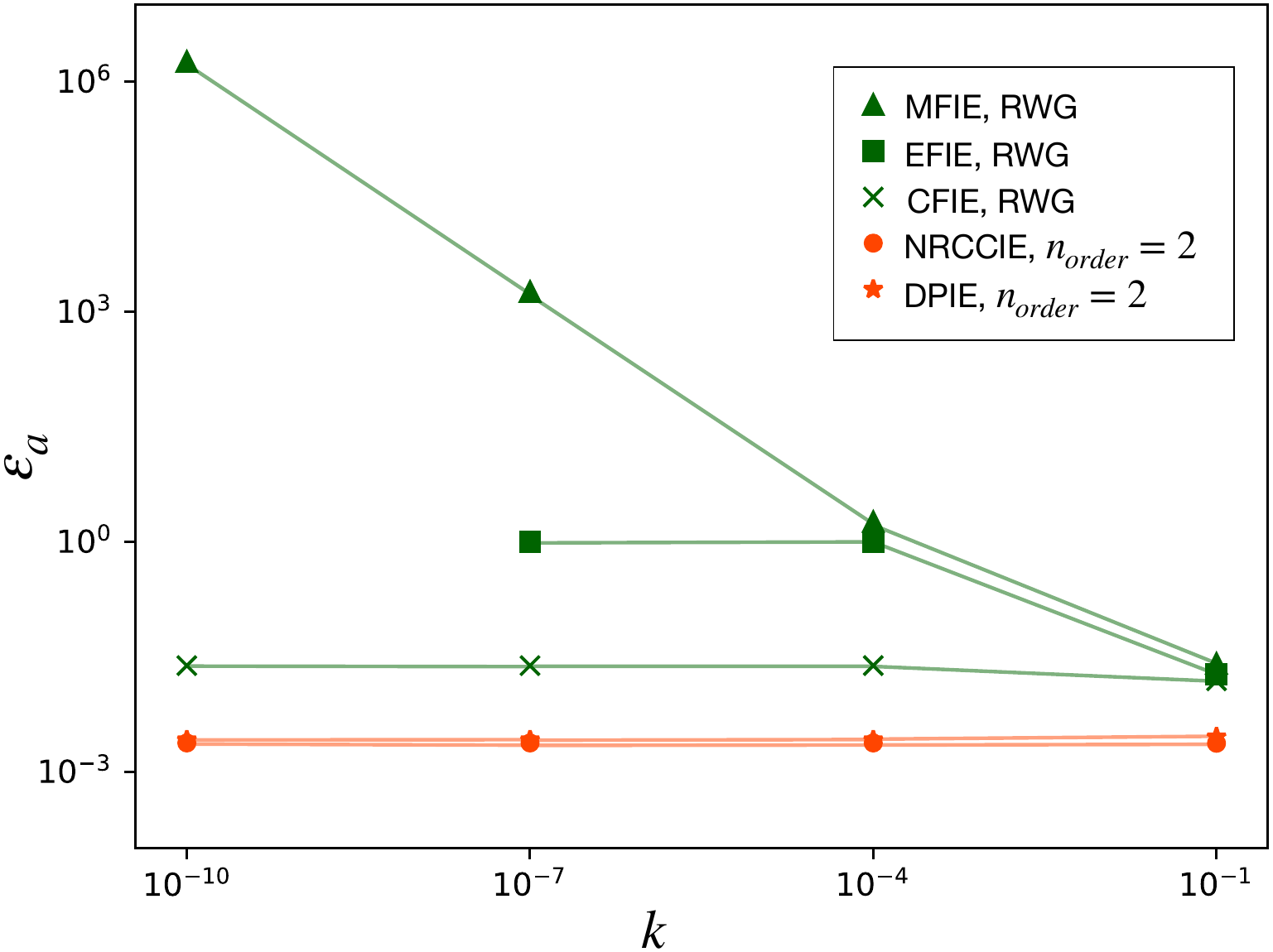}
  \caption{Relative error $\varepsilon_{a}$ as a function of wavenumber $k$ for
  the MFIE, EFIE, CFIE, NRCCIE, and DPIE on the unit sphere discretized using
  $\Ntri=192$}
  \label{fig:sphere_static-limit}
\end{center}
\end{figure}

In Fig.~\ref{fig:sphere_static-limit}, we plot the error $\varepsilon_{a}$ as
a function of $k$, with $k\in [10^{-10},10^{-1}] \SI{}{m^{-1}}$ for the MFIE, EFIE,
CFIE, NRCCIE, and DPIE. All of the integral equations were discretized using 192
patches; for the NRCCIE and DPIE we use an $\nord=2$ discretization. We
note that the CFIE, NRCCIE, and DPIE have no deterioration in accuracy in the
limit~$k\to 0$, however, for the EFIE, the error increases to~$O(1)$ as we
decrease~$k$. For the MFIE, the error increases like~$O(1/k)$ as~$k\to 0$.
\begin{remark}
Note that the MFIE is well-conditioned 
(on a surface of genus zero) in the limit $k \to 0$. However, there
is loss of accuracy in the computation of the electric field at low frequencies.
\end{remark}

\begin{figure}[t]
  \begin{center}
    \includegraphics[width=0.8\linewidth]{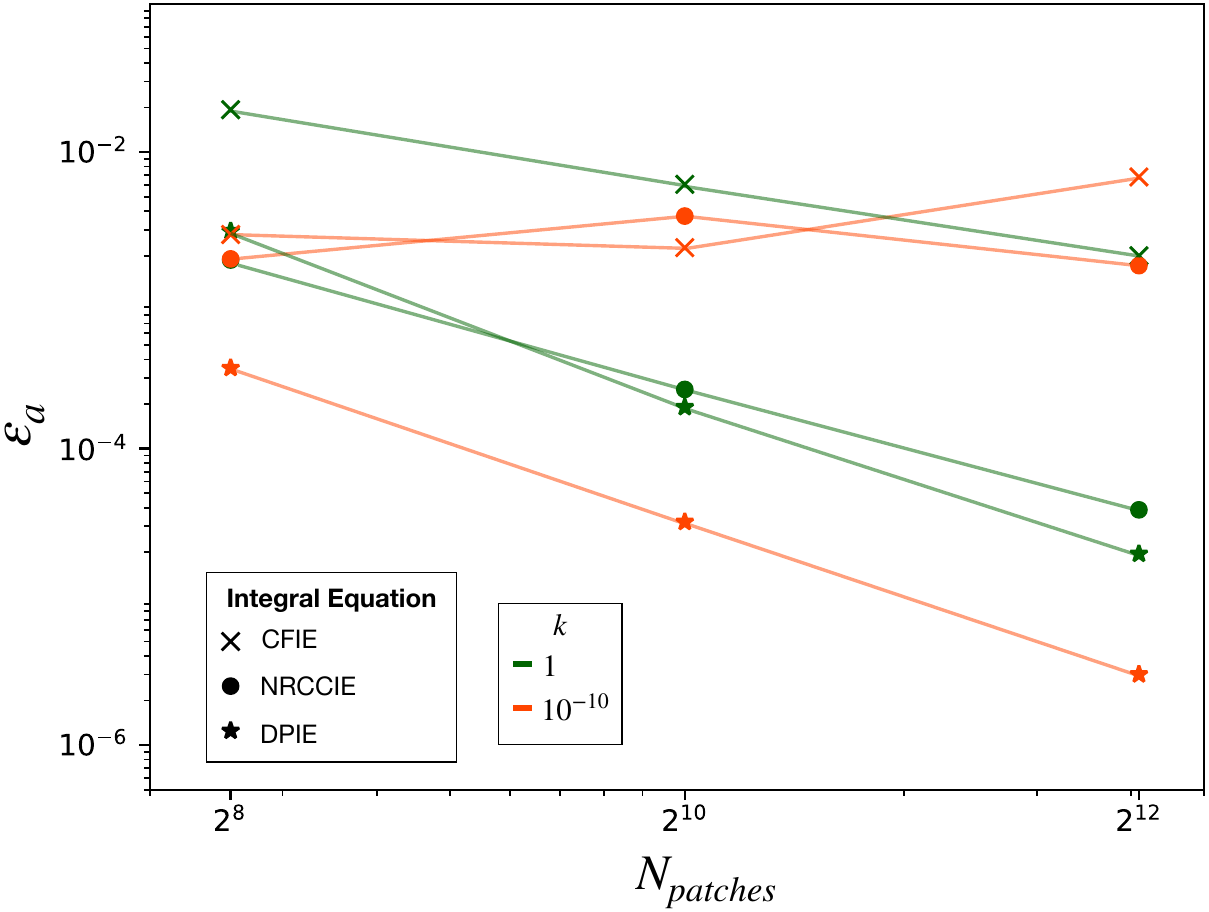}%
    \caption{Relative error $\varepsilon_{a}$ as a function of number of patches
      $\Ntri$ for the CFIE, NRCCIE, and DPIE on a smooth torus with
      $k=\SI{1}{m^{-1}}$ and $k=10^{-10}\SI{}{m^{-1}}$.}
      \label{fig:torus_static-limit}
  \end{center}
\end{figure}

The nature of the limiting static equations depends on the genus of the
conductor and the number of connected components. Thus, the stability of the
integral equation may be a function of the topology of the conductor. In
Fig.~\ref{fig:torus_static-limit}, we compare the convergence rates for the
CFIE, NRCCIE, and DPIE on the smooth torus as we refine the mesh for
$k=\SI{1}{m^{-1}}$ and $k=10^{-10}\SI{}{m^{-1}}$. The error in the computed
fields converge at the expected rate for all the integral formulations when
$k=\SI{1}{m^{-1}}$. On the other hand, for $k=10^{-10}$\SI{}{m^{-1}}, the error
in the fields computed via the DPIE continues to converge at the expected rate,
while the accuracy deteriorates upon mesh refinement for the CFIE and NRCCIE.

\begin{remark}
  \label{rem:lowfreq}
  For conductors whose dimensions are extremely small compared to the wavelength
  of the incident field, one could in principle use the solution to the static
  problems (possibly with including corrections on the Green's function~$g_{k}$
  up to~$O(k)$) in order to obtain high fidelity approximations of the
  corresponding low-frequency solutions. Such approximations are widely used in
  practice, see~\cite{haber2007octree,kapur1997ies3}, for example.
\end{remark}

\subsubsection{A posteriori error estimation}
\label{sec:post-error}

For high-order discretizations, the tail of Koornwinder expansions on each patch
can be used as an estimate for the error in the solution computed via integral
equations. Following the discussion in Section~\ref{sec:err-est},
let~$\delta J_{j}$ denote the tail of the Koornwinder expansion of the current
computed using the NRCCIE and consider the following monitor function
\begin{equation}
   \varepsilon_{j} =  \delta J_j \left( \frac{ \|\delta\bJ\|_2 / \|\bJ\|_2} {\max\limits_{j} \delta J_j} \right) .
\end{equation} 
The monitor function~$\varepsilon_{j}$ is piecewise constant on each triangle,
is proportional to~$\delta J_j$, and its maximum $\|\delta\bJ\|_2/\|\bJ\|_2$ is
the expected accuracy in the induced current. Typically, the error obtained with
the spectral monitor
function~$\varepsilon_{\rm tail} = \max_{j} \varepsilon_{j}$ is within an order
of magnitude of the relative error in the computed scattered
field~$\varepsilon_{a}$,
i.e.~$0.1 \leq \varepsilon_{a}/\varepsilon_{\rm tail} \leq 10$.

For the NRCCIE on the sphere with wavenumber $k=\SI{1}{m^{-1}}$, $\Ntri=192$,
and $\nord=4$, the estimated error from the Koornwinder tails of the current is
$\varepsilon_{\rm tail} = 1.8 \times 10^{-4}$, while the error in the field
measurements is $\varepsilon_{a} = 3.2 \times 10^{-5}$. This behavior is
independent of the wavenumber, geometry, order of discretization, number of
patches used to discretize the conductor, and also holds for other high-order
discretizations of second-kind integral equations including, the DPIE. Thus, the
error monitor function~$\varepsilon_{j}$ can reliably be used to determine
adaptive mesh refinement strategies, and~$\varepsilon_{\rm tail}$ is a
reasonable empirical indicator of the error of the solution on geometries where
an analytic solution is not known.

\subsection{Iterative solver performance}
\label{subsec:iterative}

\begin{figure}[!t]
\begin{center}
\includegraphics[width=0.8\linewidth]{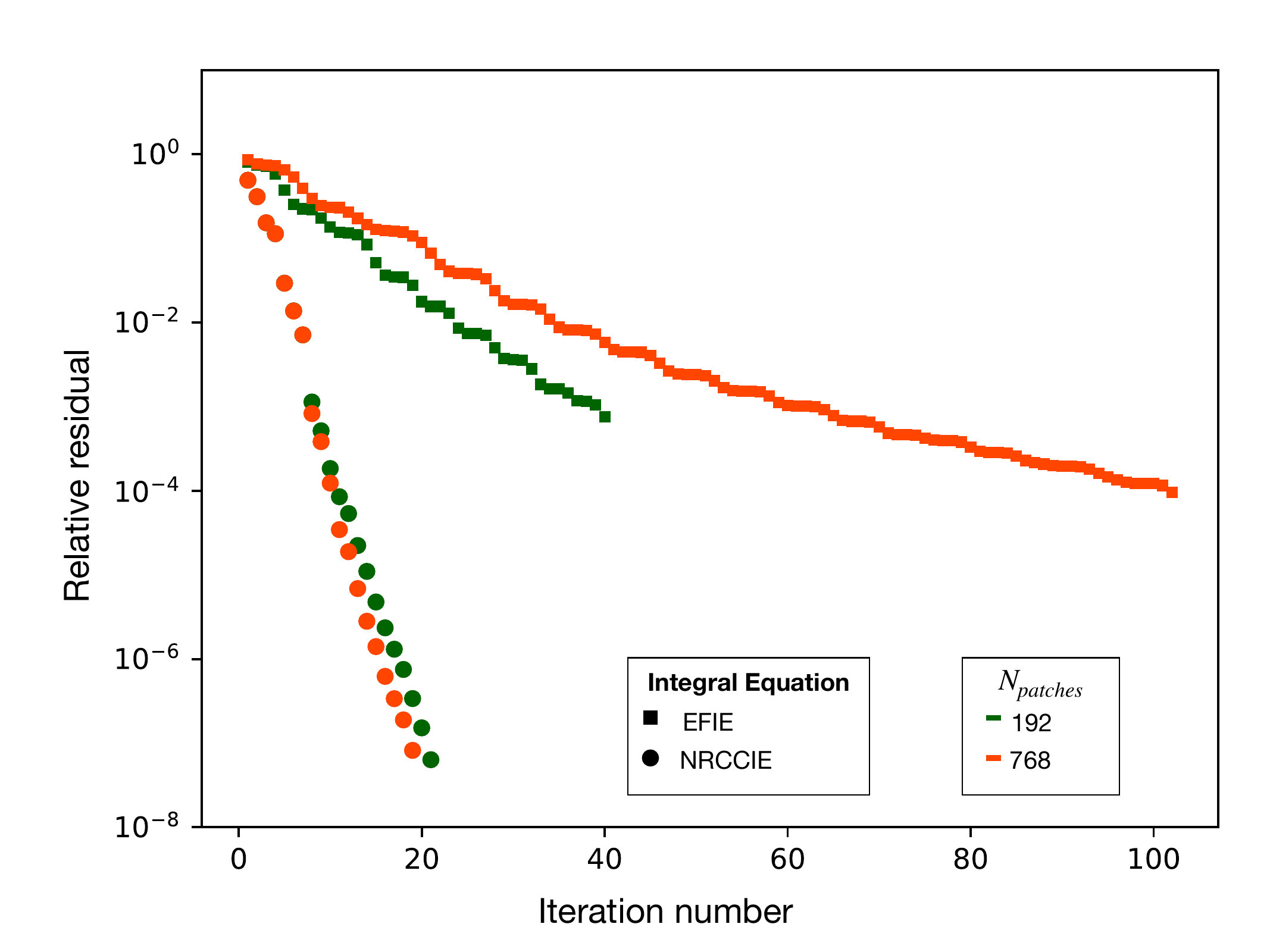}%
\end{center}
\caption{Relative GMRES residual for the NRCCIE and the EFIE.}
\label{sphere_iterations}
\end{figure}

In this section, we compare the performance of the integral equations when
coupled to an iterative solver like GMRES. It is a desirable feature for the
GMRES residual to reduce at a rate which is \emph{only} dependent on the
underlying physical problem, e.g. the complexity of the geometry and the
boundary data, but independent of the mesh used to discretize the surface. In
Fig.~\ref{sphere_iterations}, we plot the relative GMRES residual as a
function of the iteration number for the NRCCIE and EFIE. Both the integral
equations were discretized with~\mbox{$\Ntri=192$} and~\mbox{$\Ntri=768$}
patches, and second-order patches were used for discretizing the surface in the
NRCCIE. The GMRES residual as a function of iteration number is stable under
refinement for the NRCCIE, while for the EFIE, the residual decreases at a
slower rate upon mesh refinement. This stability in performance for the NRCCIE
can be attributed to its second-kind nature, while the increased number of
iterations for the EFIE can be attributed to the hypersingular nature of the
EFIE operator --- this phenomenon is often referred to as \emph{high density
mesh breakdown}. The iteration count is independent of the discretization order,
and number of patches for other second-kind integral equations, such as the DPIE
and MFIE, while integral equations with hypersingular kernels like the CFIE also
suffer from the high-density mesh breakdown.

\subsection{Large-scale examples}

We next demonstrate the performance of the integral equations on several
large-scale examples. We first demonstrate the efficiency of the DPIE on a
complicated multi-genus surface in the static limit, followed by a comparison of
the CFIE and NRCCIE for computing the far-field pattern from a bent rectangular
cavity. Finally, we illustrate the efficacy of NRCCIE for computing the
far-field pattern from a multiscale ship geometry.

\subsubsection{High-genus object in the static limit}

\begin{figure}[!t]
\begin{center}
\includegraphics[width=0.5\linewidth]{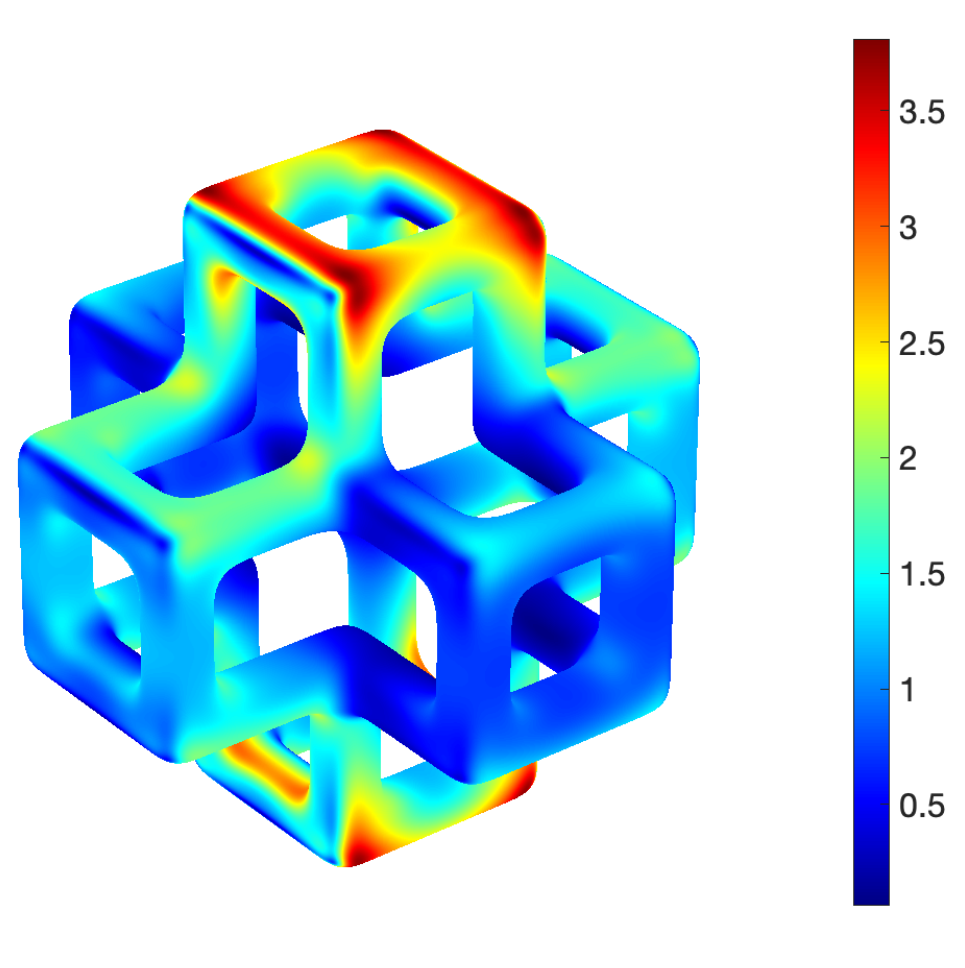}%
\end{center}
\caption{Induced source on the surface of the genus 17 geometry.}
\label{genus17_source}
\end{figure}

In the following example, we demonstrate the efficacy and stability of computing
the far-field pattern from a genus 17 surface (see Fig.~\ref{genus17_source})
in the static limit using the DPIE. None of the other integral equations
considered in this manuscript are both numerically and mathematically stable in
this regime. The incoming field is a plane wave with wavenumber
$k=10^{-10}\SI{}{m^{-1}}$. The geometry is contained in a bounding box of size
$1.6 \times 10^{-10}$ wavelengths in each dimension. As noted in
Remark~\ref{rem:lowfreq}, one could, in principle, solve a limiting PDE to
obtain a high accuracy approximation of the solution at such low frequencies.
However, computing the static solutions requires knowledge of A-cycles and
B-cycles on the geometry (i.e. loops through the holes of the surface), which
can pose a computational geometry challenge on such complicated high-order
surface meshes. The DPIE, on the other hand, can be used directly on the surface
triangulation without the need to compute these global loops on the surface.

We first compute a reference solution for this geometry where the surface is
discretized with $\nord=8$ and $\Ntri=3840$. In Fig.~\ref{genus17_source}, we
plot the induced source on the surface of the conductor using this
discretization. In order to estimate the accuracy of the computed solution and
demonstrate the efficiency of the error monitor function discussed in
Section~\ref{sec:post-error}, we also compute the solution using $\nord=2$, and
$\Ntri=960$. For this configuration, GMRES required $\niter = 72$ for the vector
equation and $\niter=21$ for the scalar part for the relative residual to reduce
to below~$10^{-6}$. The tolerance for
computing the quadrature corrections was $10^{-4}$. The accuracy in the computed
far field pattern (as measured in dB) as compared to the far field pattern
computed using the reference solution is~$1.5 \times 10^{-4}$. Another
remarkable feature of this calculation is that the DPIE can stably evaluate the
far field pattern with values ranging between $[-146,-134]$ dB, which is orders
of magnitude smaller than the induced current or the size of the conductor.

\subsubsection{Cavity}

\begin{figure}[!t]
  \centering
  \includegraphics[width=0.7\linewidth]{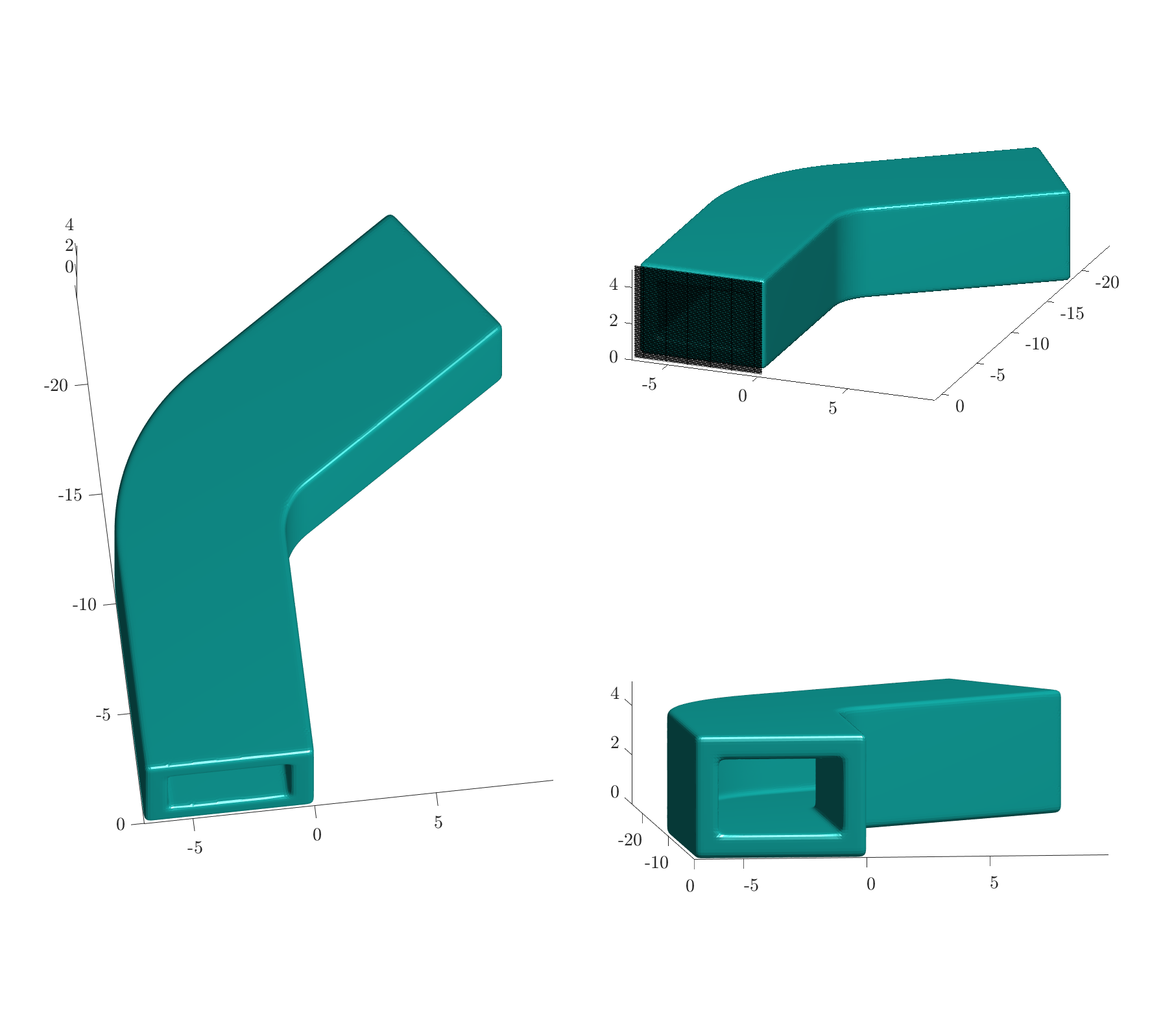}%
  \caption{Different views of the cavity. Accuracy in near-field measurements of the fields is evaluated at the screen of targets located at $x=0.5$, and $(y,z) \in [-7,0] \times[-5,0]$ as shown in the top right plot.}
  \label{Cavity_composition}
\end{figure}

\begin{figure}[t]
\begin{center}
\includegraphics[width=0.85\linewidth]{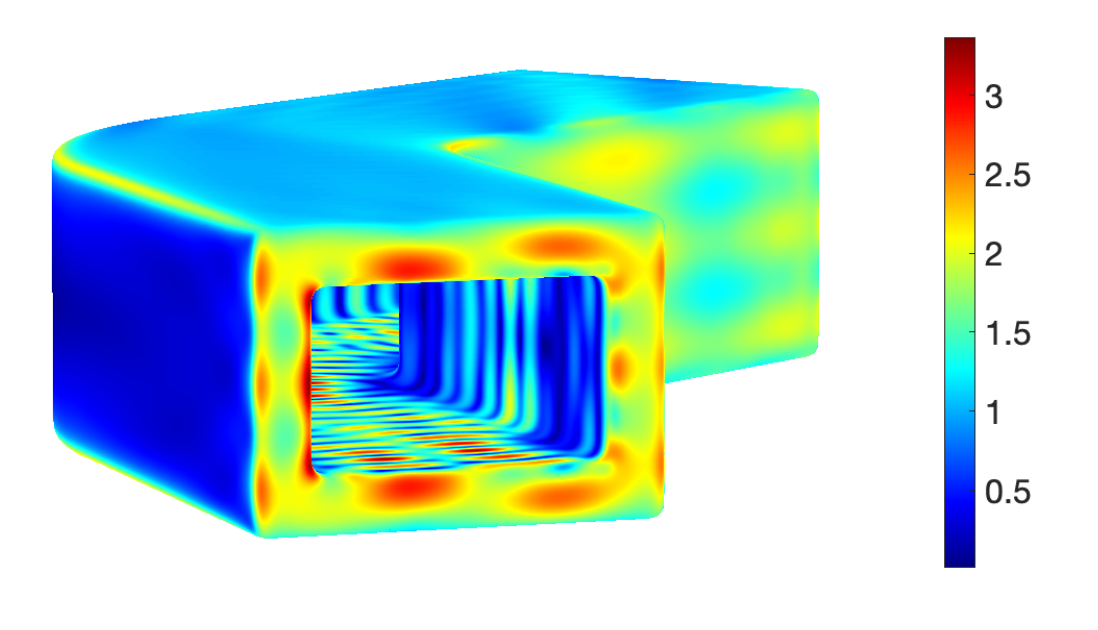}%
\end{center}
\caption{Induced current $|\bJ|$ for incoming plane wave.}
\label{Cavity_current}
\end{figure}

Next we analyze a cavity in moderately high frequency regime. The rectangular
cavity is open at one end, and around $16 \lambda$ in length along the center
line. The closed end of the cavity cannot be seen from the opening, see
Fig.~\ref{Cavity_composition}. The incoming field is a plane wave propagating
in the $-\hat{\bx}$ direction and polarized in the $\hat{\bz}$ direction. Due to
multiple internal reflections, the physical condition number of the problem is
expected to be high, and therefore this problem is a good stress test for
high-order methods. The surface of the cavity was discretized using the NRCCIE
with $(\nord,\Ntri)=(2,11392),(4,2848),(4,11392)$, and using the CFIE with
$\Ntri=11392$, and $\Ntri=45568$. The reference solution for the far field was
computed using the NRCCIE with $\Ntri=11392$ patches of $\nord=8$. The estimated
error in the reference solution based on the error monitor
function~$\varepsilon_{\textrm{tail}}=6 \times10^{-5}$. The dominant contributor
to the error in the reference solution was the tolerance used for the fast
multipole methods and quadrature corrections which was set
to~$5 \times 10^{-7}$. In Fig.~\ref{Cavity_current}, we plot the magnitude of
the induced current computed using the NRCCIE.

\begin{figure}[t]
  \begin{center}
  \includegraphics[width=\linewidth]{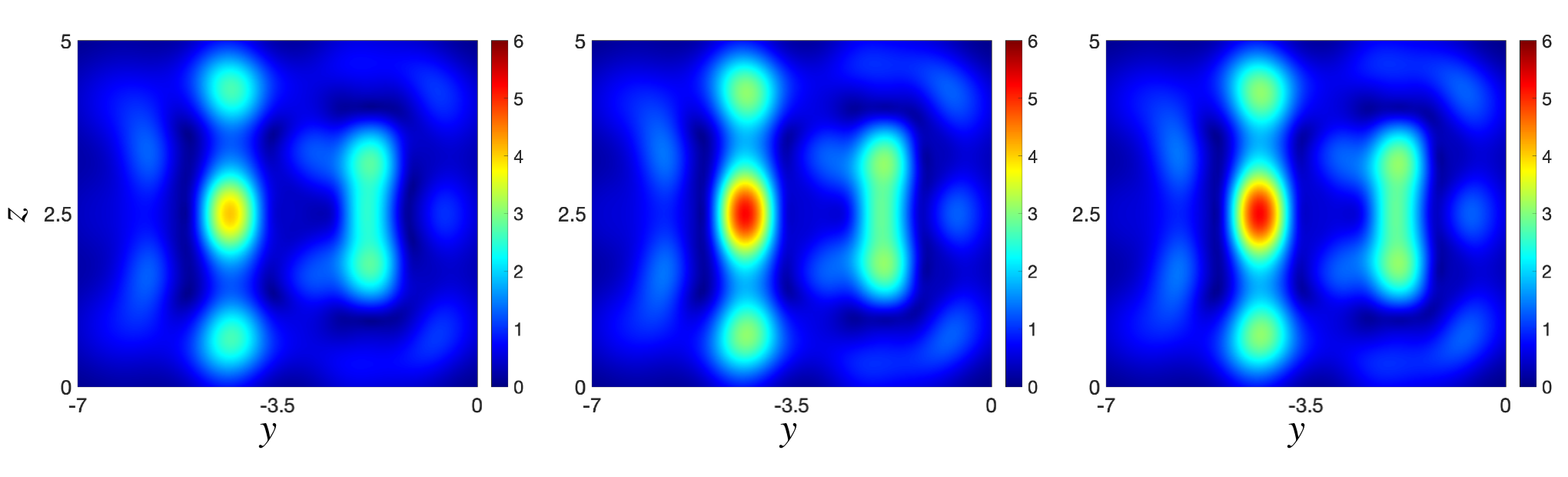}%
  \end{center}
  \caption{Norm of the Poynting vector $\|\bE^{\Sc} \times \bH^{\Sc}\|$ on the screen of targets at $x = 0.5$ and $(y,z) \in [-7,0] \times [0.5]$. (left) CFIE with $\Ntri=45568$, (middle) NRCCIE with $(\nord, \Ntri) = (4, 11392)$, and (right) reference solution.}
  \label{cavity_near}
  \end{figure}

\begin{figure}[t]
  \begin{center}
  \includegraphics[width=0.9\linewidth]{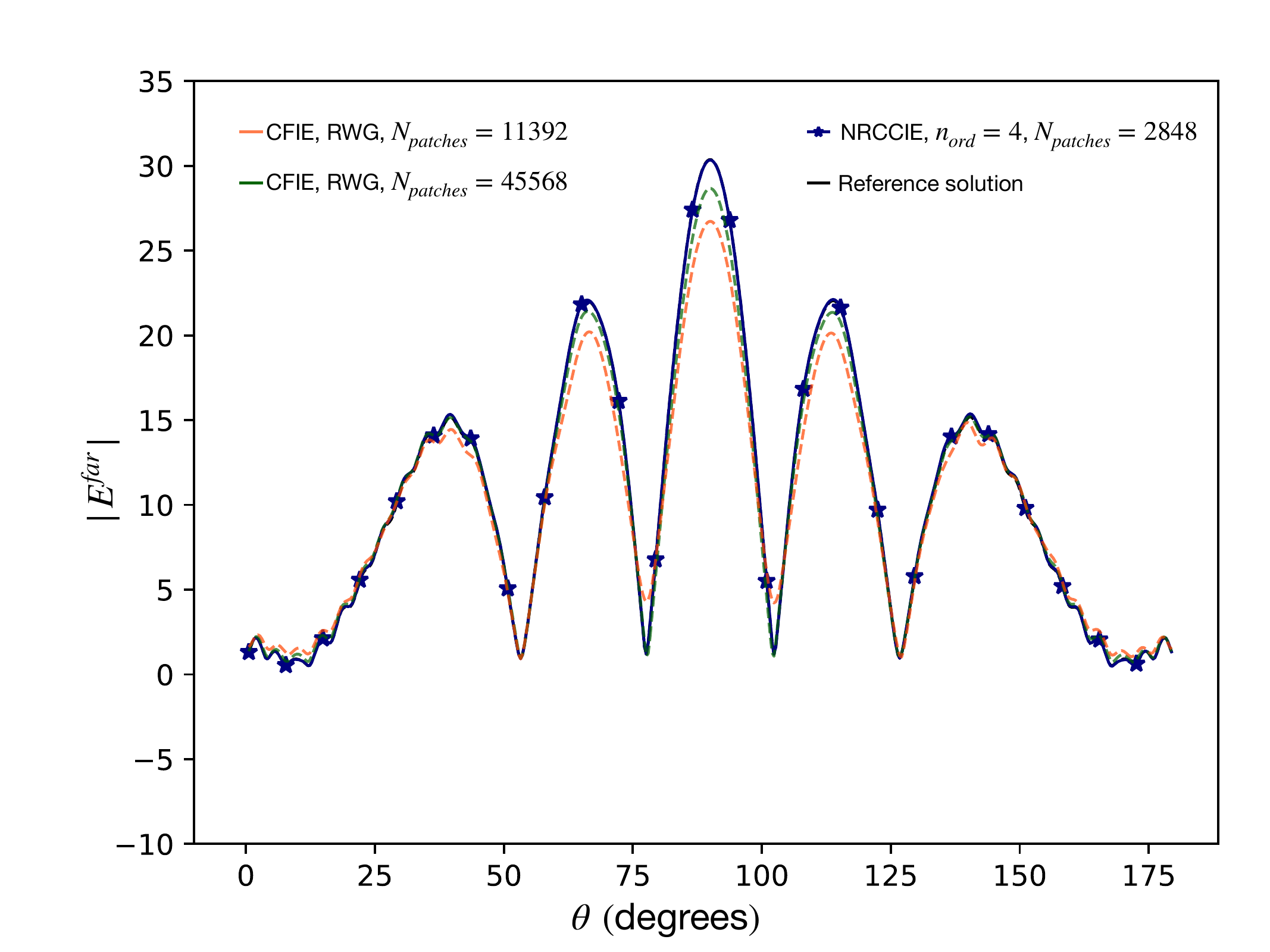}%
  \end{center}
  \caption{Far field produced by the cavity in Fig.~\ref{Cavity_composition} for
  an incoming plane wave.
  (The NRCCIE solution is indistinguishable
  in the plot from the reference solution.)}
  \label{Cavity_farfield1}
  \end{figure}

In Table~\ref{table:cavity_timing}, we tabulate the number of iterations
required for the residual to drop below the specified GMRES tolerance
$\varepsilon_{\textrm{GMRES}} = 10^{-7}$ ($\niter$ ). The precision for
computing layer potentials and the FMM was set to~$10^{-4}$. We also tabulate
the relative $L^{2}$ error in field measurements $\varepsilon_{a}$
  on a screen of targets in front of the open end of the cavity (see
  Fig.~\ref{Cavity_composition}), and the far field of the electric field
denoted by~$\varepsilon_{f}$.  In Fig.~\ref{cavity_near}, we plot
  the norm of the Poynting vector $\|\bE^{\Sc} \times \bH^{\Sc}\|$ at the near
  field targets computed using CFIE with $\Ntri = 45568$ on the left, using
  NRCCIE with $(\nord, \Ntri) = (4, 11392)$ in the middle, and the reference
  Poynting vector on the right.  In Fig.~\ref{Cavity_farfield1}, we plot the
far
field~\mbox{$\bL(\theta) = \bL\left(\sin{(\theta)},0,\cos{(\theta)} \right)$}
corresponding to the NRCCIE, the CFIE, and the reference solution for
$\theta \in [0,180]$ degrees, where $\bL(\hat{\bx})$ is as defined
in~\eqref{Fourier_NF_lf}.

\begin{table}[!b]
\begin{center}
  \caption{Iteration count, error in the near field, and error in far field
    pattern for the solution on a rectangular cavity of approximately $16$
    wavelengths along the center line due to an incoming linearly polarized
    plane wave.}
    \label{table:cavity_timing}
\begin{tabular}{c|ccccc}
& $\nord$ & $\Ntri$ & $\niter$ & $\varepsilon_{a}$ & $\varepsilon_{f}$ \\
\hline
    \multirow{3}{1.25cm}{NRCCIE} & $2$ & $11392$ & $609$  & $1 \times 10^{\scalebox{0.8}{-2}}$ & $6 \times 10^{\scalebox{0.8}{-3}}$ \\
   & $4$ & $2848$ & $607$  & $1 \times 10^{\scalebox{0.8}{-2}}$ & $5 \times 10^{\scalebox{0.8}{-3}}$ \\
  &  $4$ & $11392$ & $508$  & $2 \times 10^{\scalebox{0.8}{-3}}$ & $7 \times 10^{\scalebox{0.8}{-4}}$ \\ \hline
      \multirow{2}{1.25cm}{CFIE} & RWG & $11392$ & $216$ & $4 \times 10^{\scalebox{0.8}{-1}}$ & $1 \times 10^{\scalebox{0.8}{-1}}$ \\
   & RWG & $45568$ & $472$ & $1 \times 10^{\scalebox{0.8}{-1}}$ & $5 \times 10^{\scalebox{0.8}{-2}}$ \\
   \end{tabular}
\end{center}
\end{table}

The table highlights features of the integral equations already observed in the
previous sections with respect to the stability of the number of GMRES iterations
required for the NRCCIE, and the growth in the number of iterations required for
the CFIE. As can be seen from the plots, the error in both the near field and the far field measurements
corresponding to the CFIE with the fine mesh is still $O(1)$, while the far
field computed using the fine mesh is nearly indistinguishable from the
reference solution. 

\subsubsection{A multiscale ship simulation}
\begin{figure}[t]
  \begin{center}
  \includegraphics[width=.95\linewidth]{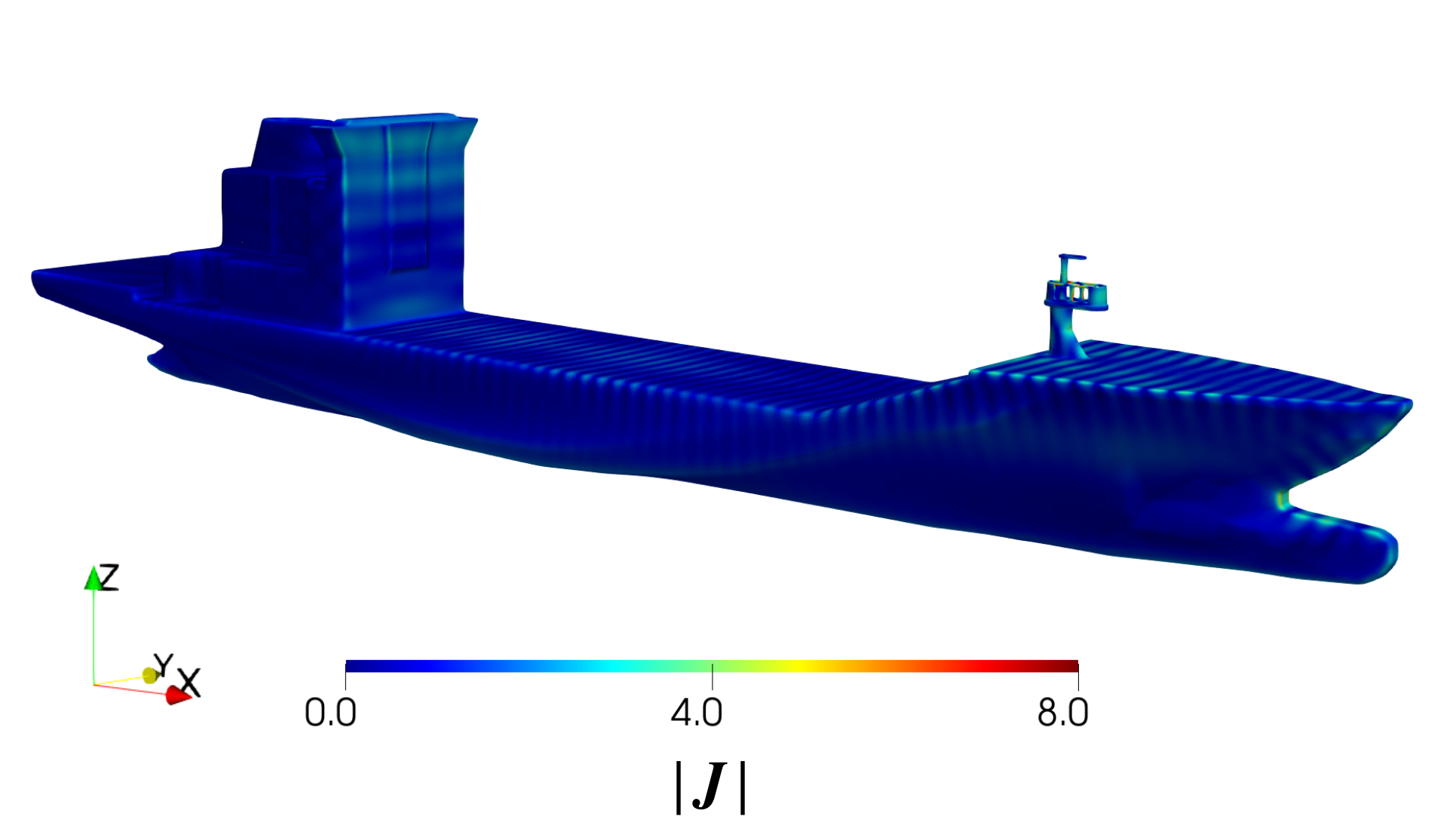}%
  \end{center}
  \caption{Induced current for incoming plane wave.}
  \label{fig:cargo_current}
\end{figure}

Finally, we demonstrate the performance of NRCCIE on a multiscale ship. The ship
is discretized using $\Ntri=30752$ patches of $\nord=4$. The incoming field is a
plane wave propagating in the~$-\hat{\bx}$ direction and polarized in the
$\hat{\bz}$ direction.  Let $R_{j}$ denote the radius of the smallest bounding
sphere containing patch $\Gamma_{j}$ centered at its centroid. The ratio of the
largest to smallest patch size, measured by the enclosing sphere radius $R_{j}$
is $21.5$. The ship is approximately $42 \lambda$ in length, $5.7 \lambda$ in
width, and $8.7 \lambda$ in height. The precision for computing the layer
potentials and the FMM were set to $10^{-4}$. For this configuration, $289$
GMRES iterations were required for the relative residual to drop below
$10^{-6}$. The estimated error in the solution is
$\varepsilon_{\textrm{tail}} = 3.8 \times 10^{-3}$. In
Fig.~\ref{fig:cargo_current}, we plot the absolute value of the induced current
on the surface of the ship.

 \section{Conclusions}

 We have demonstrated the numerical properties of various integral equation
 methods for solving exterior Maxwell scattering problems, comparing standard
 RWG discretizations of standard formulations (e.g. EFIE, MFIE, CFIE) to
 high-order Nystr\"om-like discretizations of more modern integral formulations
 especially designed to overcome the failure modes of existing ones (e.g. DPIE,
 NRCCIE).  Furthermore, we've shown that when \emph{all} aspects of the problem
 are discretized to high-order -- the geometry, quadrature, fast algorithm,
 etc. -- that high-order accuracy can be achieved at a cost which is less than
 that required for lower accuracy using standard 1st-order discretizations. We
 plan to perform a similar analysis comparing existing integral equation
 formulations with more modern ones for scattering from piecewise homogeneous
 dielectric bodies and for problems involving open surfaces (where
   the integral equation formulations are less mature).  Finally, note that we
   have omitted timing experiments in the present paper.  All implementations
   are linear scaling, but details are different for each scheme and code
   optimization involves an important but complementary set of issues.

\section*{Acknowledgments}

The authors would like to thank James Bremer, Charlie Epstein, and Zydrunas
Gimbutas for various codes and discussions during the preparation of this work. 
The Flatiron institute is a division of the Simons Foundation.


\end{document}